\documentclass[final,1p,times,twocolumn,numbers,sort&compress]{elsarticle}
\usepackage{natbib}

 \usepackage{graphicx}
\usepackage{epsfig}
\usepackage{epstopdf}

\usepackage{subfigure}
\usepackage{color}

\usepackage{latexsym,bm}

\usepackage{amssymb}
\usepackage{amsmath}
\usepackage{dsfont}

\usepackage[linesnumbered,ruled]{algorithm2e}
\usepackage{algorithmic}

\usepackage{cases}

\usepackage{threeparttable}

\newtheorem{Def}{Definition}
\newtheorem{thm}{Theorem}
\newtheorem{lem}{Lemma}
\newtheorem{cor}{Corollary}

\newcommand{\norm}[1]{\left\Vert#1\right\Vert}
\newcommand{\abs}[1]{\left| #1 \right|}
\newcommand{\pref}[1]{(\ref{#1})}

\allowdisplaybreaks

\begin{document}

\begin{frontmatter}

\title{The Convergence Rate and Necessary-and-Sufficient Condition for the Consistency of Isogeometric Collocation Method}


\author[linadd]{Hongwei Lin\corref{cor1}}
    \cortext[cor1]{Corresponding author: phone number: 86-571-87951860-8304, fax number: 86-571-88206681, email:
    hwlin@zju.edu.cn}
\address[linadd]{Department of Mathematics, State Key Lab. of CAD\&CG, Zhejiang University, Hangzhou, 310027, China}
\author[linadd]{Yunyang Xiong}

\author[huadd]{Qianqian Hu}
\address[huadd]{Department of Mathematics, Zhejiang Gongshang
                    University, Hangzhou, 310018, China}

\date{}

%\address{}
%\maketitle

\begin{abstract}
 Although the isogeometric collocation (IGA-C) method has been
    successfully utilized in practical applications due to its simplicity
    and efficiency,
    only a little theoretical results have been established on the numerical analysis of the
    IGA-C method.
 In this paper, we deduce the convergence rate of the consistency of
    the IGA-C method.
 Moreover, based on the formula of the convergence rate,
    the necessary and sufficient condition
    for the consistency of the IGA-C method is developed.
  These results advance the numerical analysis of the IGA-C method.
\end{abstract}

\end{frontmatter}

%%%%%%%%%%%%%%%%%%%%%%%%%%%%%%%%%%%%%%%%%%%%%%%%%%%%%%%%%%%%%%%%%%%%

%-----------------------------------------------------------------------------------------
% Section: Introduction
%-----------------------------------------------------------------------------------------

\section{Introduction}

 In order for the integration of CAD and CAE,
    Hughes et. al.~\cite{hughes2005isogeometric} developed the
    isogeometric analysis (IGA) method.
 Since it is based on non-linear NURBS basis
    functions,
    the IGA method can directly process the CAD models represented by NURBS,
    and avoid the tedious mesh transformation procedure.

  Because the degree of the non-linear NURBS basis function is
    relatively high, it is possible to seek a numerical solution, i.e.,
    a NURBS function, by applying the collocation method on the strong
    form of a differential equation.
  In this way, the isogeometric collocation (IGA-C) method was
    proposed~\cite{auricchio2010isogeometric}. Then unknown coefficients of the NURBS function
    can be determined by solving a linear system of equations,
    which is constructed by holding the strong form of the
    differential equation at some discrete points,
    called \emph{collocation points}.

 The IGA-C method is a simple and efficient method for solving the
    unknown coefficients of the NURBS function.
 A comprehensive
    study~\cite{schillinger2013isogeometric} revealed its superior
     behavior over the Galerkin method in terms of accuracy-to-computational-time
    ratio.
 Due to these merits, the IGA-C method has been successfully applied in some practical applications.
 However, the thorough numerical analysis for the IGA-C method is far from being established.
 Auricchio et. al. developed numerical analysis of the IGA-C
    method in one-dimensional case~\cite{auricchio2010isogeometric}.
 In the generic case, only some sufficient conditions were presented for the consistency and
    convergence of the IGA-C method~\cite{lin2013consistency}.

 In this paper, we first develop the convergence rate of the consistency of
    the IGA-C method,
    and then present the necessary and sufficient condition for the consistency of the IGA-C method.
 Specifically, for a given boundary (or initial) problem with $\mathcal{D}T=f$
    (refer to Eq.~\pref{eq:bd_p}),
    where $\mathcal{D}$ is its differential operator.
 Suppose $T_r$ is its numerical solution, represented by a NURBS function,
    and $\mathcal{I}$ is an interpolation operator such that $\mathcal{I}f = \mathcal{D}T_r$.
 The IGA-C method is consistency,
    if and only if
    %there exists a series of spline spaces,
%    in which the knot grid size series tends to $0$,
%    so that
    $\mathcal{D}$ and $\mathcal{I}$ are both uniformly bounded when the knot grid size tends to $0$.

 The structure of this paper is as follows.
 In Section~\ref{subsec:related_work},
    some related work is briefly reviewed.
 After introducing some preliminaries in Section~\ref{sec:preliminary},
    an introductory example is presented in Section~\ref{sec:intro_example}.
 Moreover, the convergence rate of the consistency of the
    IGA-C method is deduced in Section~\ref{sec:convergence_order},
    and the necessary and sufficient condition is developed in Section~\ref{sec:nsf}.
 In addition, some numerical examples are presented in
     Section~\ref{sec:numerical_examples}.
 Finally, Section~\ref{sec:conclusion} concludes the paper.

%-------------------------------------------------------------------------
% Section: related work
%-------------------------------------------------------------------------
\subsection{Related work}
\label{subsec:related_work}

 As stated above, the IGA method~\cite{hughes2005isogeometric} was proposed to
    advance the seamless integration of CAD and CAE,
    by avoiding mesh transformation.
 Moreover, since it has much less freedom than
    the traditional finite element method,
    the IGA method can not only save lots of computation,
    but also greatly improve the computational precision.
 Additionally, owing to the knot insertion property of the NURBS
    function, the shape of the computational domain represented by NURBS can be
    exactly kept in the mesh refinement.
 Due to these merits, the IGA method draws great interests in
    both practical applications and theoretical studies.
 On one hand, the IGA method has been successfully
    applied in lots of simulation problems, such as
    elasticity~\cite{auricchio2007fully,
    elguedj2008view},
    structure~\cite{cottrell2006isogeometric, hughes2008duality,
    wall2008isogeometric},
    and fluid~\cite{bazilevs2008isogeometric, bazilevs2006fluid,
    bazilevs2009patient}, etc.
 On the other hand, some research on the computational aspect of the IGA
    method has been developed to improve the accuracy and efficiency by using
    reparameterization and
    refinement, etc.~\cite{bazilevs2006isogeometric, cottrell2007studies,
    hughes2010efficient, aigner2009swept, xu2011optimal, xu2011parameterization}.
 Recently, an optimal and totally robust multi-iterative method was
    developed for solving IgA Galerkin linear system~\cite{donatelli2015robust}.
 For more details on the IGA method,
    please refer to Ref.~\cite{elguedj2014isogeometric} and the
    references therein.

 Since a NURBS function has a relatively high degree,
    its unknown coefficients can be determined by making the strong
    form of the PDE hold at some collocation points,
    that leads to the IGA-C method~\cite{auricchio2010isogeometric}.
 Schillinger et. al. presented a comprehensive comparison between
    the IGA-C method and the Galerkin method,
    revealing that the IGA-C method is superior to the Galerkin
    method in terms of accuracy-to-computational-time
    ratio~\cite{schillinger2013isogeometric}.
 Lin et. al. developed some sufficient conditions for the
    consistency and convergence of the IGA-C method~\cite{lin2013consistency}.
 Moreover, Lorenzis et. al. proposed the IGA-C method for solving
    the boundary problem with Neumann boundary condition~\cite{de2015isogeometric}.

 The IGA-C method has been successfully applied in some
    practical applications.
 For instance, the IGA-C method was employed in solving Timoshenko
    beam problem~\cite{beirao2012avoiding} and
    spatial Timoshenko rod problem~\cite{auricchio2013locking},
    showing that mixed collocation schemes are locking-free independently
    of the choice of the polynomial degrees for unknown fields.
 Moreover, the IGA-C method was extended to multi-patch NURBS configurations,
    various boundary and patch interface conditions,
    and explicit dynamic analysis~\cite{auricchio2012isogeometric}.
 Recently, the IGA-C method was exploited to settle the Bernoulli-Euler
    beam problem~\cite{reali2015isogeometric} and the Reissner-Mindlin plate
    problem~\cite{kiendl2015isogeometric}.
 However, only very limited theoretical results for the IGA-C method were
    developed~\cite{auricchio2010isogeometric, lin2013consistency} currently,
    and the numerical analysis for the IGA-C method is still far from being
    established.

%------------------------------------------------------------------------------------------
% Section: Preliminary
%------------------------------------------------------------------------------------------
\section{Preliminaries}
\label{sec:preliminary}

 Suppose the IGA-C method is employed to solve
    the following boundary problem,
 \begin{equation}\label{eq:bd_p}
    \begin{cases}
     & \mathcal{D} T = f, \qquad \qquad \text{in}\ \Omega \subset
            \mathbb{R}^d, \\
     & \mathcal{G} T = g, \qquad \qquad \text{on}\ \partial \Omega,
    \end{cases}
 \end{equation}
    where $\Omega \subset \mathbb{R}^d$ is a physical domain of $d$ dimension,
    $\mathcal{D}: \mathbb{V} \rightarrow \mathbb{W}$ is a bounded differential operator,
    where $\mathbb{V}$ and $\mathbb{W}$ are two Hilbert spaces,
    $\mathcal{G}T$ is a boundary condition,
    and $f : \Omega \rightarrow \mathbb{R}$,
    $g: \partial \Omega \rightarrow \mathbb{R}$ are two given continuous
    functions defined on their domains.
 Suppose the analytical solution $T \in C^m(\Omega)$,
    where $m$ is larger than or equal to the maximum order of derivatives appearing in the operator $\mathcal{D}$.

 In the IGA method, the physical domain $\Omega$ is represented by a
    NURBS mapping,
    \begin{equation} \label{eq:f_mapping}
        F: \Omega_p \rightarrow \Omega,
    \end{equation}
    where $\Omega_p$ is a parameter domain.
 Replacing the control points of $F$ by unknown control
    coefficients,
    we get the representation of the numerical solution
    to the boundary problem~\pref{eq:bd_p},
    denoted as $T_r(\bm{\eta}), \bm{\eta} \in \Omega_p$.
 Meanwhile, by the inverse mapping $F^{-1}$,
    the physical domain $\Omega$ can be mapped into the parameter domain
    $\Omega_p$,
    and then, the numerical solution $T_r$ is still defined on the physical domain $\Omega$ through the mapping $F^{-1}$.
 Additionally, by the mapping $F$,
    the function $f$ can be defined on $\Omega_p$,
    and $G$ on $\partial \Omega_p$.

 In isogeometric analysis, while the physical domain of the boundary
     problem~\pref{eq:bd_p} is $\Omega$,
     the computational domain is the parameter domain $\Omega_p$~\pref{eq:f_mapping}.
 Although the operators $\mathcal{D}$ and $\mathcal{G}$  in
    Eq.~\pref{eq:bd_p} are performed on the variables in the physical domain,
    the generated formulae will be transformed into the parameter domain $\Omega_p$ for computation.
 Therefore, the functions in the function approximation problem in
     the IGA-C method should be considered to be defined on the parameter domain $\Omega_p$.

%----------------------------------------------------------
% Definition:
%----------------------------------------------------------
 \begin{Def}[Stable operator~\cite{solin2006partial}]
 \label{def:stable_operator}
 Let $\mathbb{V},\ \mathbb{W}$ be Hilbert spaces and $\mathcal{D}: \mathbb{V} \rightarrow \mathbb{W}$ be a
    differential operator.
 If there exists a constant $C_S > 0$ such that
 \begin{equation*}
    \norm{\mathcal{D}v}_{\mathbb{W}} \geq C_S \norm{v}_{\mathbb{V}},\ \text{for all}\ v
    \in D(\mathcal{D}),
 \end{equation*}
 where $D(\mathcal{D})$ represents the domain of $\mathcal{D}$,
 then the differential operator $\mathcal{D}$ is called a stable
 operator.
 \end{Def}

 \textbf{[Remark 1:]} In this paper, we suppose that the $L^{\infty}$
    norm $\norm{\cdot}_{L^{\infty}}$ is equivalent to the norm
    $\norm{\cdot}_\mathbb{W}$ in $\mathbb{W}$ and the norm $\norm{\cdot}_\mathbb{V}$ in $\mathbb{V}$.
 In other words, there exists nonnegative constants $c_v, C_v$, and
    $c_w, C_w$ satisfying,
    \begin{align*}
        &c_v \norm{\cdot}_\mathbb{V} \leq \norm{\cdot}_{L^{\infty}} \leq C_v
        \norm{\cdot}_\mathbb{V} \\
        &c_w \norm{\cdot}_\mathbb{W} \leq \norm{\cdot}_{L^{\infty}} \leq C_w
        \norm{\cdot}_\mathbb{W} \\
    \end{align*}

 Suppose $T_r(\bm{\eta})$ is an unknown NURBS function defined on the knot grid
    $\mathcal{T}^{\rho} \in \mathbb{R}^d,\ d=1,2,3$.
 Specifically, $\mathcal{T}^\rho$ is a
    knot sequence in $1D$ case, a rectangular grid in $2D$ case,
    and a hexahedral grid in $3D$ case,
    where $\rho$ is the \textbf{\emph{knot grid size}} defined as the
    following definition.

 \begin{Def}\label{def:knot_grid_size_h}
    Given a set $\Phi \subset \mathbb{R}^d$, its \textbf{diameter} $diam(\Phi)$ is defined
    by
        \begin{equation*}
            diam(\Phi) = \sup\{d(\bm{x},\bm{y}), \bm{x}, \bm{y} \in \Phi\},
        \end{equation*}
        where $d(\bm{x},\bm{y})$ denotes the Euclidean distance between $\bm{x}$ and
        $\bm{y}$.
    And we call $\rho$ as the \textbf{knot grid size} of $\mathcal{T}^{\rho}$,
        which is defined as the maximum of the diameters of the \textbf{knot intervals} of $\mathcal{T}^{\rho}$.
        That is,
        $\rho = \max_i \{diam([u_i,u_{i+1}))\}$ in 1D case,
        $\rho = \max_{ij} \{diam([u_i,u_{i+1}) \times [v_j,v_{j+1}))\}$ in
        2D case,
        and $\rho = \max_{ijk} \{diam([u_i,u_{i+1}) \times [v_j,v_{j+1})
        \times [w_k, w_{k+1})) \}$ in 3D case.
 \end{Def}

 \begin{Def} \label{def:modulus_continuity}
 Let $T: \Omega_p \rightarrow \mathbb{R},\ T \in C^0(\Omega_p)$ be a
    continuous function on the parameter domain $\Omega_p$,
    where $C^0(\Omega_p)$ is the space of continuous functions on $\Omega_p$.
 The \textbf{modulus of continuity}~\cite{de2001practical} of the function $T$, denoted as $\omega(T,h)$, is defined by
 \begin{equation} \label{eq:mod_of_cont}
    \omega(T,h) = \max\{\abs{T(\bm{x})-T(\bm{y})}, d(\bm{x},\bm{y}) <
    h\}, \ h \in \mathbb{R}.
 \end{equation}
 \end{Def}

 The modulus of continuity $\omega(T,h)$ satisfies the property~\cite{de2001practical},
 \begin{equation*}
    \omega(T,h+k) \leq \omega(T,h) + \omega(T,k),\ h,k \in
    \mathbb{R},
 \end{equation*}
 and then
 \begin{equation} \label{eq:modulus}
    \omega(T,K\rho) \leq K \omega(T, \rho),\ K \in \mathbb{Z}.
 \end{equation}

 \begin{Def} \label{def:distance_of_function}
 Let $\mathcal{I}^\rho$ be an interpolation operator,
    and $\mathcal{I}^\rho g$ be a spline
    interplant of a function $g$ defined on the knot grid $\mathcal{T}^\rho$.
 %In $1D$ case, $\mathcal{I}^\rho g$ is a univariate spline function;
%    in $2D$ and $3D$ cases,
%    $\mathcal{I}^\rho g$ is a tensor product spline function.
 Suppose $\mathbb{P}$ is a spline space composed of the splines with the same knot grid
    and degree as those of $\mathcal{I}^\rho g$.
    \textbf{The distance of the function $g$ to $\mathbb{P}$}, i.e., $dist(g, \mathbb{P})$,
    is defined by
    \begin{equation} \label{eq:dist_to_space}
        dist(g,\mathbb{P}) = \min\{\norm{g-p}_{L^\infty}, p \in \mathbb{P}\}.
    \end{equation}
 \end{Def}

%--------------------------------------------------------------------------
% Section: An introductory example
%--------------------------------------------------------------------------
\section{An introductory example}
\label{sec:intro_example}

 Consider the following one-dimensional boundary problem:
 \begin{equation} \label{eq:one_dim_bd}
   \begin{cases}
       &T'(x) = f(x), \qquad \qquad \qquad \qquad x \in [a, b],\\
       &T(a) = g_1,\ T(b) = g_2,
   \end{cases}
 \end{equation}
    where $f(x) \in C[a,b]$ is a continuous function,
    $T(x) \in C^1[a,b]$ is an analytical solution,
    and $g_1,\ g_2 \in \mathbb{R}$.

 The physical domain $[a,b]$ in Eq.~\pref{eq:one_dim_bd} is modeled as,
 \begin{equation} \label{eq:one_dim_domain}
   x(t) = \sum_{i=0}^N (a + \frac{i}{N}(b-a)) B_{i,k}(t),\ t \in [0,1],
 \end{equation}
 where $B_{i,k}(t)$ is a B-spline basis function of order $k$,
 defined on the knot sequence,
 \begin{equation} \label{eq:example_knot}
   G: \underbrace{0, 0, \cdots, 0}_{k}, \frac{1}{N}, \frac{2}{N}, \cdots, \frac{N-1}{N}, \underbrace{1, 1, \cdots, 1}_{k}.
 \end{equation}
 Eq.~\pref{eq:one_dim_domain} maps $[0,1]$ to $[a,b]$, i.e.,
 \begin{equation}\label{eq:one_dim_map}
   F_1: [0,1] \rightarrow [a,b].
 \end{equation}
 Then, the numerical solution $T_r(t)$ to the boundary
     problem~\pref{eq:one_dim_bd} can be generated by replacing the coefficients $a + \frac{i}{N}(b-a)$ in $x(t)$~\pref{eq:one_dim_domain} by the unknowns coefficients $p_i\ i = 0,1,\cdots, N$, i.e.,
 \begin{equation} \label{eq:one_solution}
   T_r(t) = \sum_{i=0}^{N} p_i B_{i,k}(t).
 \end{equation}
 Note that, by the inverse mapping $F_1^{-1}$~\pref{eq:one_dim_map},
    $T_r(t)$ is defined on the physical domain $[a,b]$~\pref{eq:one_dim_bd},
    i.e., $T_r(t(x)), x \in [a,b]$.

 Because,
 \begin{equation*}
   \frac{d T_r(t)}{dt} = (k-1) \sum_{i=0}^{N-1}
            \frac{p_{i+1}-p_i}{\frac{i+k-1}{N}-\frac{i}{N}}
            B_{i,k-1}(t)
            =
            N \sum_{i=0}^{N-1} (p_{i+1}-p_i) B_{i,k-1}(t),
 \end{equation*}
 and,
 \begin{equation*}
   \frac{d x(t)}{d t} = b-a,
 \end{equation*}
    substituting Eq.~\pref{eq:one_solution} into Eq.~\pref{eq:one_dim_bd} yields,
 \begin{equation*}
   \begin{cases}
    & \frac{d T_r}{dx} = \frac{d T_r}{dt} \frac{dt}{dx}
        = \frac{d T_r}{dt} \frac{1}{\frac{dx}{dt}}
        =  \sum_{i=0}^{N-1} \frac{N}{b-a} (p_{i+1}-p_i) B_{i,k-1}(t)
        = f(x(t)), \\
    &T_r(0) = p_0 = g_1,\ T_r(1) = p_{N} = g_2.
   \end{cases}
 \end{equation*}

 In order for solving the unknown coefficients in Eq.~\pref{eq:one_solution}
    using the IGA-C method,
    a linear system is generated by sampling $N-1$ points
    $\tau_1, \tau_2, \cdots, \tau_{N-1}$ in
    the interval $(0,1)$, i.e.,
 \begin{equation} \label{eq:one_d_linear_system}
   \begin{cases}
    & \frac{dT_r(\tau_j)}{dx}=
        \sum_{i=0}^{N-1} \frac{N}{b-a} (p_{i+1}-p_i) B_{i,k-1}(\tau_j)
        = f(x(\tau_j)),
    \ \tau_j \in (0,1),\ j = 1,2,\cdots,N-1, \\
    &T_r(0) = p_0 = g_1,\\
    &T_r(1) = p_n = g_2.
   \end{cases}
 \end{equation}
 When the knot grid size $\rho = \frac{1}{N}$ of the knot sequence
     $G$~\pref{eq:example_knot} tends to $0$,
    it follows $N \rightarrow +\infty$.
 If the control points $\frac{N}{b-a}(p_{i+1}-p_i) \rightarrow \infty,
    (\rho \rightarrow 0)$, too,
    we have $\frac{dT_r(\tau_j)}{dx} \rightarrow \infty, (\rho \rightarrow 0), j = 1, 2, \cdots, N-1$.
 However, because $f(x(t))$ is continuous on the close interval $[0,1]$,
    $f(x(\tau_j))$ is bounded.
 Therefore, if the linear system~\pref{eq:one_d_linear_system} has
    a solution,
    there should exist $T_r(t)$ so that the control points
    $\frac{N}{b-a}(p_{i+1}-p_i)$ of $\frac{dT_r}{dx}$ are bounded when $\rho = \frac{1}{N} \rightarrow 0$.
 It results in that $\frac{dT_r}{dx}$ is also bounded when
    $\rho \rightarrow 0$.
 All of such B-spline functions $T_r(t)$ constitute a B-spline subspace,
    and the first order derivative operator in Eq.~\pref{eq:one_dim_bd}
    should be bounded on the B-spline subspace when $\rho \rightarrow 0$.

%--------------------------------------------------------------------------
% Section: The convergence rate
%--------------------------------------------------------------------------
\section{The convergence rate}
\label{sec:convergence_order}

 Suppose the NURBS function
    $T_r(\bm{\eta}),\ \bm{\eta} \in \Omega_p \subset \mathbb{R}^d$
    defined on the knot grid $\mathcal{T}^\rho$
    has $n$ unknown control coefficients $p_{\bm{i}}$, i.e.,
    \begin{equation}\label{eq:rewrite-nurbs}
    T_r(\bm{\eta}) = \sum_{\bm{i}} p_{\bm{i}} \frac{w_{\bm{i}} B_{\bm{i}}(\bm{\eta})}{W(\bm{\eta})}
                   = \frac{P(\bm{\eta})}{W(\bm{\eta})},
                    \quad \bm{\eta} \in \Omega_p \subset \mathbb{R}^d,
 \end{equation}
    where $w_{\bm{i}} > 0$ are known weights,
    $B_{\bm{i}}(\bm{\eta})$ are the B-spline basis functions,
    the weight function $W(\bm{\eta})$ is a known polynomial spline function,
    and $P(\bm{\eta})$ is a polynomial spline function with
    $n$ unknown control coefficients $p_{\bm{i}}$.
 Moreover, the subscript $\bm{i}$ in Eq.~\pref{eq:rewrite-nurbs} is an index vector, $\bm{i} = (i_1, i_2, \cdots,
    i_d)$.
 According to the IGA-C method, these unknown coefficients $p_{\bm{i}}$ can be determined by solving the
    following linear system of equations,
    \begin{equation}\label{eq:linear_system}
        \begin{cases}
        & \mathcal{D}T_r (\bm{\eta}_k) = f(\bm{\eta}_k),\qquad
            \qquad k = 1,2,\cdots,n_1, \\
        & \mathcal{G}T_r (\bm{\eta}_l) = g(\bm{\eta}_l), \qquad
            \qquad\ l = n_1 + 1, \cdots, n,
        \end{cases}
    \end{equation}
    where $\bm{\eta}_k (k = 1,2,\cdots,n_1)$ are collocation points
    inside $\Omega_p$,
    and $\bm{\eta}_l (l=n_1+1, \cdots, n)$ are collocation points on
    $\partial \Omega_p$.
 Note that, throughout this paper, the operators $\mathcal{D}$ and
    $\mathcal{G}$ are performed on the variable in the physical domain $\Omega$ (Eq.~\pref{eq:bd_p}).

 \vspace{0.3cm}

 \textbf{[Remark 2:]} In this paper, we assume that the coefficient matrix of the above linear
    system~\pref{eq:linear_system}
    is of full rank and then it has a unique solution.
 Otherwise, the IGA-C method is invalid.

 \vspace{0.3cm}

 According to the result developed in Ref.~\cite{lin2013consistency},
    $\mathcal{D}T_r$ can be represented by
    \begin{equation} \label{eq:dtr_representation}
        \mathcal{D}T_r(\bm{\eta})
            = \sum_{\bm{i}} p_{\bm{i}} \mathcal{D} \frac{w_{\bm{i}}B_{\bm{i}}(\bm{\eta})}{W(\bm{\eta})}
            = \sum_{\bm{i}} p_{\bm{i}} \frac{\bar{B}_{\bm{i}}(\bm{\eta})}{\bar{W}(\bm{\eta})}
            = \frac{\bar{P}(\bm{\eta})}{\bar{W}(\bm{\eta})},
    \end{equation}
 where $\bar{B}_i(\bm{\eta})$ is the result by applying the
    differential operator $\mathcal{D}$ to
    $\frac{w_{\bm{i}}B_{\bm{i}}(\bm{\eta})}{W(\bm{\eta})}$,
    $\bar{W}(\bm{\eta})$ is the power of $W(\bm{\eta})$,
    and $\bar{P}(\bm{\eta})$ is a polynomial B-spline function with $n$ unknowns $p_{\bm{i}}$.
 By Ref.~\cite{lin2013consistency}, $\bar{P}(\bm{\eta})$ and $\bar{W}(\bm{\eta})$ both have the same
    break point sequence and the same knot intervals as $T_r(\bm{\eta})$.
 To determine these unknowns $p_{\bm{i}}$ in $\bar{P}(\bm{\eta})$,
    let $\mathcal{D} T_r (\bm{\eta})$
    interpolate $\mathcal{D}T(\bm{\eta}) = f(\bm{\eta})$ at $n_1$
    collocation points inside the domain $\Omega_p$ (refer to Eq.~\pref{eq:linear_system}),
    i.e.,
    \begin{equation} \label{eq:op_equivalence}
        \mathcal{D}T_r(\bm{\eta}_k) - f(\bm{\eta}_k) =
        \frac{\bar{P}(\bm{\eta}_k)}{\bar{W}(\bm{\eta}_k)} - f(\bm{\eta}_k)
        = \frac{\bar{P}(\bm{\eta}_k) -
        \bar{W}(\bm{\eta}_k)f(\bm{\eta}_k)}{\bar{W}(\bm{\eta}_k)} = 0,
        \ k = 1, 2, \cdots, n_1.
    \end{equation}
 Note that $\bar{W}(\bm{\eta}) \neq 0$ is a known function,
    Eq.~\pref{eq:op_equivalence} is equivalent to
    \begin{equation} \label{eq:dtr_pbar}
        \bar{P}(\bm{\eta}_k) = \sum_{\bm{i}} p_{\bm{i}}\bar{B}_{\bm{i}}(\bm{\eta}_k)
        = \bar{W}(\bm{\eta}_k)f(\bm{\eta}_k),\ k = 0,1,\cdots,n_1.
    \end{equation}

 Similarly, $\mathcal{G}T_r (\bm{\eta})$ in
    Eq.~\pref{eq:linear_system} can be written as
    \begin{equation} \label{eq:gtr_representation}
        \mathcal{G}T_r (\bm{\eta})
        = \sum_{\bm{i}} p_{\bm{i}} \mathcal{G}
            \frac{w_{\bm{i}}B_{\bm{i}}(\bm{\eta})}{W(\bm{\eta})}
        = \sum_{\bm{i}} p_{\bm{i}} \frac{\tilde{B}_{\bm{i}}(\bm{\eta})}{\tilde{W}(\bm{\eta})}
        = \frac{\tilde{P}(\bm{\eta})}{\tilde{W}(\bm{\eta})},
    \end{equation}
    where $\tilde{B}_{\bm{i}}(\bm{\eta})$ are the result generated by
    applying the operator $\mathcal{G}$ to
    $\frac{w_{\bm{i}}B_{\bm{i}}(\bm{\eta})}{W(\bm{\eta})}$,
    $\tilde{W}(\bm{\eta}) \neq 0$ is a known B-spline
    function,
    and $\tilde{P}(\bm{\eta})$ is an unknown B-spline
    function with $n$ unknowns $p_{\bm{i}}$.
 Then the linear equations $\mathcal{G}T_r (\bm{\eta}_l) = g(\bm{\eta}_l)$
    in Eq.~\pref{eq:linear_system} are equivalent to
    \begin{equation} \label{eq:gtr_ptilde}
        \tilde{P}(\bm{\eta}_l) = \sum_{\bm{i}} p_{\bm{i}}\tilde{B}_{\bm{i}}(\bm{\eta}_l)
        = \tilde{W}(\bm{\eta}_l) g(\bm{\eta}_l),\ l = n_1+1, n_1+2, \cdots, n,
    \end{equation}
    where $\bm{\eta}_l \in \partial \Omega_p,\ l=n_1+1, \cdots, n$.

 Therefore, combining Eqs.~\pref{eq:dtr_pbar} and~\pref{eq:gtr_ptilde},
    the linear system~\pref{eq:linear_system} becomes
  \begin{equation}\label{eq:equ_linear_system}
        \begin{cases}
        & \bar{P}(\bm{\eta}_k)
            = \sum_{\bm{i}} p_{\bm{i}}\bar{B}_{\bm{i}}(\bm{\eta}_k)
            = \bar{W}(\bm{\eta}_k)f(\bm{\eta}_k),\ k = 0,1,\cdots,n_1, \\
        & \tilde{P}(\bm{\eta}_l)
            = \sum_{\bm{i}} p_{\bm{i}}\tilde{B}_{\bm{i}}(\bm{\eta}_l)
            = \tilde{W}(\bm{\eta}_l)
            g(\bm{\eta}_l),\ l = n_1+1, n_1+2, \cdots, n.
        \end{cases}
    \end{equation}
  Since the linear system of equations~\pref{eq:equ_linear_system} is
    equivalent to~\pref{eq:linear_system},
    then the coefficient matrix of~\pref{eq:equ_linear_system} is of full
    rank, and it also has a unique solution.

\vspace{0.3cm}

 \textbf{[Remark 3:]} In Eq.~\pref{eq:linear_system}, let the functions $f$ and $g$ vary
    in $C^m(\Omega_p)$ and $C^m(\partial \Omega_p)$, respectively,
    and the differential operator $\mathcal{D}$ be fixed.
 In addition, let the weight function $W(\bm{\eta})$ in $T_r$~\pref{eq:rewrite-nurbs}
    be fixed as well.
 Then, all the numerical solutions
    $T_r(\bm{\eta})$~\pref{eq:rewrite-nurbs}
    generated by the IGA-C method~\pref{eq:linear_system}
    constitute \textbf{a linear spline space} $\mathbb{S}_{\rho}(\Omega_p)$,
    where $\rho$ is the knot grid size of $T_r$.
 It should be pointed out that,
    all the NURBS functions $T_r$ in the linear space $\mathbb{S}_{\rho}(\Omega_p)$
    have the same weight function $W(\bm{\eta})$,
    the same knot grid with knot grid size $\rho$ and the same degree.
 In order for $\rho \rightarrow 0$,
    the knot grid of the spline functions in $\mathbb{S}_{\rho}(\Omega_p)$ is refined by knot insertion,
    thus resulting in a series of spline spaces.
 %Moreover, $\mathcal{D} T_r$ can be represented as Eq.~\pref{eq:dtr_representation},
%    i.e.,
%    $\mathcal{D}T_r(\bm{\eta}) =
%    \frac{\bar{P}(\bm{\eta})}{\bar{W}(\bm{\eta})}$.
 Moreover, because all the numerical solutions $T_r(\bm{\eta})$ constitute
    the linear space $\mathbb{S}_{\rho}(\Omega_p)$,
    all of $\mathcal{D}T_r(\bm{\eta})$ compose \textbf{a linear spline space}
    $\mathbb{S}^d_{\rho,e}(\Omega_p)$,
    where $e$ is the continuity order of the splines in $\mathbb{S}^d_{\rho,e}(\Omega_p)$.
 As aforementioned,
    $\mathcal{D}T_r$ has the same break point sequence with that of
    $T_r$~\pref{eq:rewrite-nurbs},
    so they have the same knot grid size $\rho$.

 \vspace{0.3cm}

  The following Lemma~\ref{lem:d_m_relation} estimates the distance from a
    continuous function $f \in C^0(\Omega_p)$ to the linear space $\mathbb{S}^d_{\rho,e}(\Omega_p)$,
    i.e., $dist(f, \mathbb{S}^d_{\rho,e})$.
 In Ref.~\cite[pp.146]{de2001practical},
    an inequality to estimate the distance is
    proposed for univariate functions,
    and the inequality can be extended to our case.

 %---------------------------------------------------------------------------
 % Lemma:
 %---------------------------------------------------------------------------
 \begin{lem} \label{lem:d_m_relation}
 If $\mathcal{D}T = f \in C^0(\Omega_p)$, and $T \in C^m(\Omega_p)$
    (Eq.~\pref{eq:bd_p}), then we have
 \begin{equation*}
    dist(f, \mathbb{S}^d_{\rho,e}) = dist(\mathcal{D}T, \mathbb{S}^d_{\rho,e}) \leq \norm{\mathcal{D}} K \omega(T,\rho),
 \end{equation*}
 where $K$ is an integer related to the degree of the NURBS functions
 in the spline space $\mathbb{S}_{\rho}(\Omega_p)$.
 \end{lem}

 \textbf{Proof:}
 As stated above, the NURBS functions approximating the analytical solution
    $T$ constitute the linear space $\mathbb{S}_{\rho}(\Omega_p)$ defined on the knot grid $\mathcal{T}^{\rho}$.
 We select a special function from the space $\mathbb{S}_{\rho}(\Omega_p)$,
    i.e.,
    \begin{equation*}
        T_r(\bm{\eta}) = \sum_{\bm{i}} T(\bm{\tau}_{\bm{i}})
            \frac{w_{\bm{i}}B_{\bm{i}}(\bm{\eta})}{W(\bm{\eta})},
    \end{equation*}
 %Referring to Eq.~\pref{eq:rewrite-nurbs},
%    suppose the NURBS function $T_r \in \mathbb{S}_{\rho}(\Omega_p)$ defined on the knot grid $\mathcal{T}^{\rho}$ and expressed as
%    \begin{equation*}
%        T_r(\bm{\eta}) = \sum_{\bm{i}} p_{\bm{i}}
%            \frac{w_{\bm{i}}B_{\bm{i}}(\bm{\eta})}{W(\bm{\eta})},
%    \end{equation*}
% approximates the analytical solution $T$.
    and construct a spline function $(Af)(\bm{\eta})$ to approximate the function
    $f \in C^0(\Omega_p)$, i.e.,
    \begin{equation*}
        (Af)(\bm{\eta}) = \mathcal{D}T_r(\bm{\eta}) = \mathcal{D} \sum_{\bm{i}}
        T(\bm{\tau}_{\bm{i}})
        \frac{w_{\bm{i}}B_{\bm{i}}(\bm{\eta})}{W(\bm{\eta})}
        = \sum_{\bm{i}}
        T(\bm{\tau}_{\bm{i}}) \mathcal{D}
        \frac{w_{\bm{i}}B_{\bm{i}}(\bm{\eta})}{W(\bm{\eta})},
    \end{equation*}
    where $T$ is the analytical solution of Eq.~\pref{eq:bd_p},
    and $Af = \mathcal{D}T_r \in \mathbb{S}^d_{\rho,e}(\Omega_p)$ (defined in Remark 3).
 The point sequence $\{\bm{\tau}_{\bm{i}} \in \Omega_p\}$ is sampled in such
    a way that each knot interval of the knot grid $\mathcal{T}^{\rho}$
    contains at least one point,
    and $\bm{\tau}_{\bm{i}}$ is in the non-zero region of $B_{\bm{i}}(\bm{\eta})$.

 Suppose $u(\bm{\eta}),\ v(\bm{\eta}) \in C^m(\Omega_p)$.
 The function $u(\bm{\eta})$ is an arbitrary function in $C^m(\Omega_p)$,
    and,
    \begin{equation} \label{eq:v_function}
      v(\bm{\eta}) = T(\bm{\eta}) - \sum_{\bm{i}} T(\bm{\tau}_{\bm{i}})
        \frac{w_{\bm{i}} B_{\bm{i}}(\bm{\eta})}{W(\bm{\eta})},
        \ \bm{\eta} \in \Omega_p.
    \end{equation}
 Note that $\abs{v(\bm{\eta})}$ is continuous in the close set $\Omega_p$,
    so $\abs{v(\bm{\eta})}$ can take its maximum value in $\Omega_p$.
 Namely, there exists $\bm{\eta}^{\ast} \in \Omega_p$ such that
 \begin{equation*}
   \abs{v(\bm{\eta}^{\ast})} = \max_{\bm{\eta}}\abs{v(\bm{\eta})} = \norm{v}_{L^\infty},\ \bm{\eta} \in \Omega_p.
 \end{equation*}
 For an arbitrary value $\hat{\bm{\eta}} \in \Omega_p$,
    it holds,
 \begin{equation*}
   \norm{\mathcal{D}} = \max_{u \in C^m(\Omega_p)} \frac{\norm{\mathcal{D}u}_{L^{\infty}}}
                                    {\norm{u}_{L^{\infty}}}
                      \geq \frac{\norm{\mathcal{D}v}_{L^\infty}}
                                {\norm{v}_{L^\infty}}
                      \geq \frac{\abs{\mathcal{D}v(\hat{\bm{\eta}})}}
                                {\abs{v(\bm{\eta}^{\ast})}},
 \end{equation*}
 which is equivalent to,
 \begin{equation} \label{eq:d_inequality}
   \abs{\mathcal{D}v(\hat{\bm{\eta}})} \leq \norm{\mathcal{D}} \abs{v(\bm{\eta}^{\ast})}.
 \end{equation}

 Suppose $J$ is the index vector set satisfying $B_{\bm{i}}(\bm{\eta}^{\ast}) \neq 0,\ \bm{i} \in J$.
 Because
 \begin{equation*}
   \sum_{\bm{i}} \frac{w_{\bm{i}}
        B_{\bm{i}}(\bm{\eta}^{\ast})}{W(\bm{\eta}^{\ast})}
 = \sum_{\bm{i} \in J} \frac{w_{\bm{i}}
        B_{\bm{i}}(\bm{\eta}^{\ast})}{W(\bm{\eta}^{\ast})}
 = 1,\
 \text{and then}\ T(\bm{\eta}^{\ast}) = \sum_{\bm{i}} T(\bm{\eta}^{\ast}) \frac{w_{\bm{i}}
        B_{\bm{i}}(\bm{\eta}^{\ast})}{W(\bm{\eta}^{\ast})}
        = \sum_{\bm{i} \in J} T(\bm{\eta}^{\ast}) \frac{w_{\bm{i}}
        B_{\bm{i}}(\bm{\eta}^{\ast})}{W(\bm{\eta}^{\ast})},
 \end{equation*}
 together with Eq.~\pref{eq:d_inequality}, we have,
 \begin{align*}
    \abs{f(\hat{\bm{\eta}}) - (Af)(\hat{\bm{\eta}})}
    &=
        \abs{\mathcal{D}T(\hat{\bm{\eta}}) - \mathcal{D}
        \sum_{\bm{i}} T(\bm{\tau}_{\bm{i}}) \frac{w_{\bm{i}}
        B_{\bm{i}}(\hat{\bm{\eta}})}{W(\hat{\bm{\eta}})}}\\
    &= \abs{\mathcal{D}v(\hat{\bm{\eta}})}
    \leq \norm{\mathcal{D}} \abs{v(\bm{\eta}^{\ast})}
    \qquad \text{(Eq.~\pref{eq:d_inequality})}\\
    & =
        \norm{\mathcal{D}} \abs{T(\bm{\eta}^{\ast}) - \sum_{\bm{i}} T(\bm{\tau}_{\bm{i}}) \frac{w_{\bm{i}}
        B_{\bm{i}}(\bm{\eta}^{\ast})}{W(\bm{\eta}^{\ast})}}
        \qquad \text{(Eq.~\pref{eq:v_function})}\\
    & =
        \norm{\mathcal{D}} \abs{\sum_{\bm{i}} T(\bm{\eta}^{\ast}) \frac{w_{\bm{i}}
        B_{\bm{i}}(\bm{\eta}^{\ast})}{W(\bm{\eta}^{\ast})} - \sum_{\bm{i}} T(\bm{\tau}_{\bm{i}}) \frac{w_{\bm{i}}
        B_{\bm{i}}(\bm{\eta}^{\ast})}{W(\bm{\eta}^{\ast})}}\\
    & =
        \norm{\mathcal{D}} \abs{\sum_{\bm{i} \in J} (T(\bm{\eta}^{\ast}) - T(\bm{\tau}_i)) \frac{w_{\bm{i}}
        B_{\bm{i}}(\bm{\eta}^{\ast})}{W(\bm{\eta}^{\ast})} } \\
    & \leq
        \norm{\mathcal{D}} \sum_{\bm{i} \in J}
        \abs{T(\bm{\eta}^{\ast})-T(\bm{\tau}_{\bm{i}})} \frac{w_{\bm{i}}
        B_{\bm{i}}(\bm{\eta}^{\ast})}{W(\bm{\eta}^{\ast})} \\
    & \leq
        \norm{\mathcal{D}} \max_{\bm{i} \in J}
        \abs{T(\bm{\eta}^{\ast})-T(\bm{\tau}_{\bm{i}})} \qquad (\text{Definition}~\ref{def:modulus_continuity})\\
    & \leq
        \norm{\mathcal{D}} \omega(T, K\rho), %\\
 \end{align*}
 where $K$ is an integer related to the degree of
 $T_r(\bm{\eta})$.
 It is because that the non-zero region of $B_{\bm{i}}(\bm{\eta})$ is determined by the degree of $T_r(\bm{\eta})$.
 By Eq.~\pref{eq:modulus}, we get
 \begin{equation*}
 \abs{f(\hat{\bm{\eta}}) - (Af)(\hat{\bm{\eta}})} \leq
        \norm{\mathcal{D}} K \omega(T, \rho).
 \end{equation*}
 Because $Af = \mathcal{D}T_r \in \mathbb{S}^d_{\rho,e}(\Omega_p)$,
    and $\hat{\bm{\eta}} \in \Omega_p$ is an arbitrary value,
    it can be so chosen that
 \begin{equation*}
    dist(f,\mathbb{S}^d_{\rho,e}) = min\{\norm{f-s}_{L^{\infty}}, s \in \mathbb{S}^d_{\rho,e}\}
    \leq \abs{f(\hat{\bm{\eta}}) - (Af)(\hat{\bm{\eta}})}
    \leq \norm{\mathcal{D}} K \omega(T, \rho).
 \end{equation*}
 Then the Lemma is proved.
 $\Box$

 \vspace{0.3cm}

 Furthermore, we have
 \begin{lem} \label{lem:estimation_m}
  If $T(\bm{\eta}) \in C^1(\Omega_p)$, then it holds
    \begin{equation*}
    \omega(T,\rho) \leq \rho \max_{\bm{\eta} \in \Omega_p} \norm{\nabla T}_E,
    \end{equation*}
    where $\nabla T$ is the gradient of $T$,
    and the norm $\norm{\cdot}_E$ is defined as $\norm{\bm{\eta}}_E = \norm{(\eta_1,\eta_2,\cdots,\eta_d)}_E = \sqrt{\eta_1^2+\eta_2^2+\cdots+\eta_d^2}$.
 \end{lem}

 \textbf{Proof:} Let $\bm{x}, \bm{y} \in \Omega_p$,
    $d(\bm{x},\bm{y}) = \norm{\bm{x}-\bm{y}}_E \leq \rho$, and $c \in (0,1)$.
 According to the mean value theorem, it follows that
 \begin{equation*}
 \begin{split}
    \abs{T(\bm{x}) - T(\bm{y})} & = \abs{\nabla T |_{(1-c)\bm{x} +
    c\bm{y}} \cdot (\bm{x}-\bm{y})}
    \leq \norm{\nabla T |_{(1-c)\bm{x} +
    c\bm{y}}}_{E} \cdot \norm{\bm{x}-\bm{y}}_{E} \\
    & \leq \rho \norm{\nabla T |_{(1-c)\bm{x} +
    c\bm{y}}}_{E}
    \leq \rho \max_{\bm{\eta} \in \Omega_p} \norm{\nabla T}_E.
 \end{split}
 \end{equation*}
 Then, by the definition of $\omega(T,\rho)$
    (Eq.~\pref{eq:mod_of_cont}), we have
 \begin{equation*}
    \omega(T,\rho) = \max\{\abs{T(\bm{x})-T(\bm{y})}, d(\bm{x},\bm{y}) < \rho \}
    \leq \rho \max_{\bm{\eta} \in \Omega_p} \norm{\nabla T}_E.
 \end{equation*} $\Box$

\vspace{0.2cm}

 Moreover, we denote by $\mathcal{I}^\rho$ an interpolation
    operator,
    which maps a continuous function to a spline function defined on
    the knot grid $\mathcal{T}^\rho$ with knot grid size $\rho$.
 Specifically, for the continuous function $f = \mathcal{D}T \in C^0(\Omega_p)$(refer
    to Eqs.~\pref{eq:bd_p} and~\pref{eq:linear_system}),
    we have
    \begin{equation} \label{eq:inter_operator}
        \mathcal{I}^\rho f = \mathcal{D}T_r \in \mathbb{S}^d_{\rho,e}(\Omega_p),
    \end{equation}
 and the following Lemma.

 %-------------------------------------------------------------
 % Theorem
 %-------------------------------------------------------------
 \begin{lem} \label{lem:bound}
    Suppose $\mathcal{D}T = f \in C^0(\Omega_p)$ (Eq.~\pref{eq:bd_p}),
        and $T_r \in \mathbb{S}_{\rho}(\Omega_p)$ (Refer to Eq.~\pref{eq:rewrite-nurbs} and Remark 3)
        is the NURBS function approximating the analytical solution $T$.
    Then,
    \begin{equation}\label{eq:}
        \norm{\mathcal{D}T_r - \mathcal{D}T}_\mathbb{W}
        \leq
        (1 + \norm{\mathcal{I}^\rho}) dist(f,\mathbb{S}^d_{\rho,e}),
    \end{equation}
    where $\mathcal{I}^\rho$ is the interpolation operator defined
    by Eq.~\pref{eq:inter_operator},
    and $\mathbb{S}^d_{\rho,e}$ is defined as in Remark 3.
 \end{lem}

 \textbf{Proof:} On one hand,
 given an arbitrary known NURBS function $T_q(\bm{\eta}) \in \mathbb{S}_{\rho}(\Omega_p)$ expressed as
 \begin{equation} \label{eq:known_fun_sufficiency}
    T_q(\bm{\eta}) = \frac{Q(\bm{\eta})}{W(\bm{\eta})}
    = \sum_{\bm{i}} q_{\bm{i}} \frac{w_{\bm{i}}
    B_{\bm{i}}(\bm{\eta})}{W(\bm{\eta})}, \bm{\eta} \in \Omega_p,
 \end{equation}
    where the weight function $W(\bm{\eta})$ and the weight $w_i$ are the same as those in
    \pref{eq:rewrite-nurbs},
    two functions $h(\bm{\eta})$ and $h_b(\bm{\eta})$ can be
    generated by performing the operators
    $\mathcal{D}$ and $\mathcal{G}$ on $T_q$ (see
    Eq.~\pref{eq:bd_p}), respectively, i.e.,
 \begin{equation} \label{eq:fun_h}
    h(\bm{\eta}) = \mathcal{D}T_q(\bm{\eta})
    = \sum_{\bm{i}} q_{\bm{i}} \mathcal{D} \left( \frac{w_{\bm{i}}
        B_{\bm{i}}(\bm{\eta})}{W(\bm{\eta})} \right),
    \
    \ h_b(\bm{\eta}) = \mathcal{G}T_q(\bm{\eta})
    = \sum_{\bm{i}} q_{\bm{i}} \mathcal{G} \left( \frac{w_{\bm{i}}
        B_{\bm{i}}(\bm{\eta})}{W(\bm{\eta})} \right).
 \end{equation}

 We construct an unknown NURBS function $T_x(\bm{\eta}) \in \mathbb{S}_{\rho}(\Omega_p)$ with $n$
    unknown control coefficients $x_{\bm{i}}$,
    the same knot grid and degree with $T_q$,
 \begin{equation} \label{eq:unknow_fun_sufficiency}
    T_x(\bm{\eta}) = \frac{X(\bm{\eta})}{W(\bm{\eta})}
    = \sum_{\bm{i}} x_{\bm{i}} \frac{w_{\bm{i}}
    B_{\bm{i}}(\bm{\eta})}{W(\bm{\eta})},
 \end{equation}
 where the weight function $W(\bm{\eta})$ and the weight $w_i$ are the same as those in
    Eqs.~\pref{eq:rewrite-nurbs} and~\pref{eq:known_fun_sufficiency}.
 The $n$ unknown coefficients $x_{\bm{i}}$ in
    $T_x(\bm{\eta})$ can be obtained by making $\mathcal{D}T_x$ and $\mathcal{G}T_x$ interpolate
    $h(\bm{\eta})$ and $h_b(\bm{\eta})$ at some sampling points, respectively,
    similar as~\pref{eq:linear_system}, i.e.,
    \begin{equation} \label{eq:system_suffi}
        \begin{cases}
        & \mathcal{D}T_x (\bm{\eta}_k) = h(\bm{\eta}_k),\qquad
            \qquad k = 1,2,\cdots,n_1, \\
        & \mathcal{G}T_x (\bm{\eta}_l) = h_b(\bm{\eta}_l), \qquad
            \qquad l = n_1 + 1, \cdots, n.
        \end{cases}
    \end{equation}
 Therefore, $\mathcal{I}^\rho h = \mathcal{D}T_x$.

 The aforementioned linear system of equations~\pref{eq:system_suffi} can be rewritten
    as
    \begin{equation}\label{eq:trans_system_suffi}
        \begin{cases}
        & \sum_{\bm{i}} (x_{\bm{i}} - q_{\bm{i}}) \mathcal{D}\left. \left( \frac{w_{\bm{i}}
        B_{\bm{i}}(\bm{\eta})}{W(\bm{\eta})} \right)\right|_{\bm{\eta} = \bm{\eta}_k} = 0,\qquad
            \qquad k = 1,2,\cdots,n_1, \\
        & \sum_{\bm{i}} (x_{\bm{i}} - q_{\bm{i}}) \mathcal{G}\left. \left( \frac{w_{\bm{i}}
        B_{\bm{i}}(\bm{\eta})}{W(\bm{\eta})} \right)\right|_{\bm{\eta} = \bm{\eta}_l} = 0, \qquad
            \qquad l = n_1 + 1, \cdots, n.
        \end{cases}
    \end{equation}
 Obviously, the coefficient matrix of~\pref{eq:trans_system_suffi} is
    the same as that of the linear system \pref{eq:linear_system},
    and is of full rank, too.
 Then the linear system of equations~\pref{eq:trans_system_suffi} has only
    zero solution, i.e., $x_{\bm{i}} =
     q_{\bm{i}}$, meaning that
 \begin{equation} \label{eq:h_in_sd}
    \mathcal{I}^\rho h = \mathcal{D}T_x = \mathcal{D}T_q = h.
 \end{equation}

 %On the other hand, note that $\mathbb{S}^d_{\rho,e}$ (refer to Remark 3) is the linear space composed of
%    the NURBS functions with the same knot grid $\mathcal{T}^{\rho}$ and degree
%    as those of $\mathcal{D}T_r$~\pref{eq:dtr_representation},
%    and
%    \begin{equation}\label{eq:fun_h}
%        h(\bm{\eta}) = \mathcal{D}T_q(\bm{\eta})
%            \in \mathbb{S}^d_{\rho,e}.
%    \end{equation}
 Therefore, we have
 \begin{align} \label{eq:consistency}
    \norm{\mathcal{D}T_r - \mathcal{D}T}_\mathbb{W} & =
    \norm{\mathcal{I}^\rho f - f}_\mathbb{W} =
    \norm{\mathcal{I}^\rho f - \mathcal{I}^\rho h + h - f}_\mathbb{W}
    = \norm{\mathcal{I}^\rho (f-h) - (f-h)}_\mathbb{W} \\
    & \leq \norm{f-h}_\mathbb{W} + \norm{\mathcal{I}^\rho} \norm{f-h}_\mathbb{W}
    = (1 + \norm{\mathcal{I}^\rho}) \norm{f-h}_\mathbb{W} \notag \\
    & = (1 + \norm{\mathcal{I}^\rho}) dist(f,\mathbb{S}^d_{\rho,e}). \quad \text{(explained in the following paragraph)} \notag
 \end{align}

  Because $T_q(\bm{\eta})$~\pref{eq:known_fun_sufficiency} is an
    arbitrary NURBS function in the spline space $\mathbb{S}_{\rho}(\Omega_p)$,
    the function $h(\bm{\eta}) = \mathcal{D}T_q$~\pref{eq:fun_h} is also an arbitrary NURBS function in the linear space $\mathbb{S}^d_{\rho,e}(\Omega_p)$.
 So in Eq.\pref{eq:consistency},
    the function $h$ can be chosen from $\mathbb{S}^d_{\rho,e}(\Omega_p)$ to make $\norm{f-h}_\mathbb{W}$ as small as possible,
    that is, $dist(f,\mathbb{S}^d_{\rho,e})$. $\Box$

 \vspace{0.3cm}

 Based on Lemma~\ref{lem:d_m_relation} and~\ref{lem:bound}, it follows:
 %-------------------------------------------------------------
 % Theorem
 %-------------------------------------------------------------
 \begin{lem} \label{lem:consistency}
    Suppose $\mathcal{D}T = f \in C^0(\Omega_p)$ (Eq.~\pref{eq:bd_p}),
        and $T_r \in \mathbb{S}_{\rho}(\Omega_p)$ (Refer to Eq.~\pref{eq:rewrite-nurbs} and Remark 3)
        is the NURBS function approximating the analytical solution $T$.
    Then,
    \begin{equation*}
        \norm{\mathcal{D}T_r - \mathcal{D}T}_\mathbb{W}
        \leq K\norm{\mathcal{D}}(1 + \norm{\mathcal{I}^\rho}) \omega(T, \rho)  ,
    \end{equation*}
    where $\mathcal{I}^\rho$ is the interpolation operator defined
    by Eq.~\pref{eq:inter_operator},
    and $K$ is an integer related to the degree of the splines in the spline space $\mathbb{S}_{\rho}(\Omega_p)$.
 \end{lem}

\vspace{0.3cm}

 Moreover, due to Lemma~\ref{lem:estimation_m} and~\ref{lem:consistency},
    the convergence rate of $\mathcal{D}T_r$ to $\mathcal{D}T$
    when $\rho \rightarrow 0$ is obtained as follows.
 \begin{thm} \label{thm:rate_consistency}
    Suppose the analytical solution $T \in C^1(\Omega_p)$ (Eq.~\pref{eq:bd_p}).
    We have,
    \begin{equation*}
        \norm{\mathcal{D}T_r - \mathcal{D}T}_\mathbb{W} \leq
        \rho K {\norm{\mathcal{D}}} (1 + \norm{\mathcal{I}^\rho})
        \max_{\bm{\eta} \in \Omega_p} \norm{\nabla T}_E.
    \end{equation*}
 Here, $\mathcal{D}T, T_r, \mathcal{I}^\rho$, and $K$ are delineated as in Lemma~\ref{lem:consistency}.
 \end{thm}

 In addition, if $\mathcal{D}$ is a stable operator (Definition~\ref{def:stable_operator}),
    we can get the convergence rate of $T_r$ to $T$ when $\rho
    \rightarrow 0$.
 \begin{cor} \label{thm:rate_convergence}
 Suppose the operator $\mathcal{D}$ in Eq.~\pref{eq:bd_p} is a stable
    differential operator,
    and $T \in C^1(\Omega_p)$.
 We have,
    \begin{equation*}
        \norm{T_r - T}_\mathbb{V}
        \leq
        \frac{K}{C_S} \norm{\mathcal{D}} (1 + \norm{\mathcal{I}^{\rho}})
        \omega(T, \rho)
        \leq
        \frac{\rho K}{C_S} \norm{\mathcal{D}} (1 + \norm{\mathcal{I}^{\rho}})
        \max_{\bm{\eta} \in \Omega_p} \norm{\nabla T}_E,
    \end{equation*}
    where $C_S$ is a positive constant, $T_r, \mathcal{I}^\rho$, and $K$ are delineated as in Lemma~\ref{lem:consistency}.
 \end{cor}

 %-----------------------------------------------------------------------------
 % \Subsection
 %-----------------------------------------------------------------------------
 \subsection{One dimensional case}
 \label{subsec:one_dim}

 In the one dimensional case,
    the convergence rate can be improved.
 In this section, suppose the operator $\mathcal{D}$ is a linear differential
     operator with constant coefficients.

 %----------------------------------------------------------------------------
 % Lemma
 %----------------------------------------------------------------------------
 \begin{lem}\cite[pp. 148]{de2001practical}
 \label{lem:dist_relation}
  Let $g \in C^m(\Omega_p)$ be a univariate function,
    and $\mathbb{S}^d_{\rho,e}(\Omega_p)$ be defined as in Remark 3.
  It holds,
  \begin{equation}\label{eq:dist_ineq}
    dist(g, \mathbb{S}^d_{\rho,e}) \leq \gamma \rho\ dist(g', \mathbb{S}^d_{\rho, e-1}),
  \end{equation}
  where $\gamma$ is a number related to the degree of the splines in $\mathbb{S}^d_{\rho,e}(\Omega_p)$,
  and $g'$ is the first order derivative of $g$.
 \end{lem}

 Repeatedly using Lemma~\ref{lem:dist_relation} leads to:
 %--------------------------------------------------------------
 % Lemma
 %--------------------------------------------------------------
 \begin{lem} \label{lem:one_dim_dist}
 Suppose $f = \mathcal{D}T \in C^m(\Omega_p)$ (Eq.~\pref{eq:bd_p}) is a
    univariate function,
    the linear spline space $\mathbb{S}^d_{\rho, e}(\Omega_p)$ is defined as in Remark 3, and the operator $\mathcal{D}$ is a linear differential operator with constant coefficients.
 We have,
 \begin{equation*}
   dist(f, \mathbb{S}^d_{\rho, e})
   = dist(\mathcal{D}T, \mathbb{S}^d_{\rho,e})
   \leq \Gamma \norm{\mathcal{D}} \rho^{\nu} \norm{T^{(\nu)}}_{L^{\infty}},
 \end{equation*}
 where $\nu = min(m,e)$, $\Gamma$ is a number related to $\nu$ and the degree
 of the splines in $\mathbb{S}^d_{\rho, e}(\Omega_p)$,
 and $T^{(\nu)}$ is the $\nu^{th}$ order derivative of $T$.
 \end{lem}

 \textbf{Proof:}
 Because $f = \mathcal{D}T \in C^m(\Omega_p)$ is a univariate function,
    and $\mathcal{D}$ is a linear differential operator with constant coefficients,
    we have $(\mathcal{D} T)^{(k)} = \mathcal{D} T^{(k)},\ k = 1,2,\cdots, m$.
 By using Lemma~\ref{lem:dist_relation} repeatedly,
    and denoting $\nu = min(m,e)$,
    it follows,
 \begin{equation*}
   \begin{split}
    dist(f, \mathbb{S}^d_{\rho,e})
       & = dist(\mathcal{D}T, \mathbb{S}^d_{\rho,e})
        \leq \gamma_1 \rho\ dist((\mathcal{D}T)', \mathbb{S}^d_{\rho,e-1})
        \leq \gamma_1 \gamma_2 \rho^2\ dist((\mathcal{D}T)'', \mathbb{S}^d_{\rho,e-2}) \\
       & \leq \cdots
        \leq  \gamma_1 \gamma_2 \cdots \gamma_{\nu-1}
            \rho^{\nu-1}\ dist((\mathcal{D}T)^{(\nu-1)}, \mathbb{S}^d_{\rho,e-\nu+1}) \\
       & = \gamma_1 \gamma_2 \cdots \gamma_{\nu-1}
            \rho^{\nu-1}\ dist(\mathcal{D}T^{(\nu-1)}, \mathbb{S}^d_{\rho,e-\nu+1})  \\
       & \leq \gamma_1 \cdots \gamma_{\nu-1} \rho^{\nu-1} K_{\nu}
            \norm{\mathcal{D}} \omega(T^{(\nu-1)},\rho), \qquad
            (Lemma~\ref{lem:d_m_relation})
   \end{split}
 \end{equation*}
 where, $\gamma_1$ is a number related to the degree of the splines in
    $\mathbb{S}^d_{\rho,e}$ (denoted as $deg$),
    $\gamma_2$ is a number related to the degree of the splines in
    $\mathbb{S}^d_{\rho,e-1}$, i.e., $deg - 1$, $\cdots$, and so on;
    $K_{\nu}$ is a number related to the degree of the splines in
    $\mathbb{S}^d_{\rho,e-\nu+1}$, i.e., $deg - \nu + 1$.
 In conclusion, $\gamma_i, i = 1,2, \cdots, \nu-1$, and $K_{\nu}$
    are all related to $deg$ and $\nu = min(m,e)$,
    and then we denote $\Gamma = \gamma_1 \gamma_2 \cdots \gamma_{\nu-1} K_{\nu}$.
 Moreover, by Lemma~\ref{lem:estimation_m}, we have,
 \begin{equation*}
   dist(f, \mathbb{S}^d_{\rho,e})  \leq \Gamma \rho^{\nu-1}
            \norm{\mathcal{D}} \omega(T^{(\nu-1)},\rho)
            \leq \Gamma \norm{\mathcal{D}} \rho^{\nu} \norm{T^{(\nu)}}_{L^\infty},
 \end{equation*}
 where, $\nu = min(m,e)$, and $\Gamma$ is a number related to $\nu$ and the degree of the splines in $\mathbb{S}^d_{\rho, e}(\Omega_p)$.
 $\Box$

 \vspace{0.3cm}

 Based on Lemma~\ref{lem:bound} and~\ref{lem:one_dim_dist},
    the convergence rate for the consistency of the IGA-C method in the one-dimensional case is deduced.
 \begin{thm} \label{thm:one_dim_consistency}
 Suppose $f = \mathcal{D}T \in C^m(\Omega_p)$ (Eq.~\pref{eq:bd_p}) is a
    univariate function,
    the spline space $\mathbb{S}^d_{\rho, e}(\Omega_p)$ is defined as in Remark 3,
    and the operator $\mathcal{D}$ is a linear differential operator with constant coefficients.
 We have,
 \begin{equation*}
   \norm{\mathcal{D}T_r - \mathcal{D}T}_{\mathbb{W}}
   \leq \Gamma (1 + \norm{\mathcal{I}^\rho}) \norm{\mathcal{D}} \rho^{\nu} \norm{T^{(\nu)}}_{L^{\infty}},
 \end{equation*}
 where $\nu = min(m,e)$, and $\Gamma$ is a number related to $\nu$ and the degree of the splines in $\mathbb{S}^d_{\rho, e}(\Omega_p)$.
 \end{thm}

 Moreover, if the operator $\mathcal{D}$ is also a stable operator
    (Definition 1),
    it holds:
 \begin{cor} \label{thm:one_dim_convergence}
 Suppose $f = \mathcal{D}T \in C^m(\Omega_p)$ (Eq.~\pref{eq:bd_p}) is a
    univariate function,
    the spline space $\mathbb{S}^d_{\rho, e}(\Omega)$ is defined as in Remark 3,
    and the linear differential operator with constant coefficients $\mathcal{D}$ is stable (refer to Definition 1).
 We have,
 \begin{equation*}
   \norm{T_r - T}_{\mathbb{V}}
   \leq \frac{\Gamma}{C_S}(1 + \norm{\mathcal{I}^\rho}) \norm{\mathcal{D}} \rho^{\nu} \norm{T^{(\nu)}}_{L^{\infty}},
 \end{equation*}
 where $C_S$ is a positive constant,
    $\nu = min(m,e)$,
    and $\Gamma$ is a number related to $\nu$ and the degree of the splines in $\mathbb{S}^d_{\rho, e}(\Omega_p)$.
 \end{cor}

%-----------------------------------------------------------------------------
% \Section
%-----------------------------------------------------------------------------
 \section{The necessary and sufficient condition}
 \label{sec:nsf}

 In this section, we will present the necessary and sufficient condition of
    the consistency of the IGA-C method.
 Because $\mathcal{D}T = f \in C^0(\Omega_p)$ and $T$ is continuous
    (Eq.~\pref{eq:bd_p}),
    we have $\omega(T,\rho) \rightarrow 0$,
    when $\rho \rightarrow 0$.
 Based on Lemma~\ref{lem:consistency},
    if $\norm{\mathcal{I}^\rho}$ and $\norm{\mathcal{D}}$ are bounded,
    it follows
    $\norm{\mathcal{D}T_r - \mathcal{D}T}_\mathbb{W} \rightarrow 0$
    when $\rho \rightarrow 0$.
 That is, the IGA-C method is consistency.
 However, since $\mathcal{I}^{\rho}f = \mathcal{D}T_r \in
    \mathbb{S}^d_{\rho,e}$~\pref{eq:inter_operator},
    and $T_r \in \mathbb{S}_{\rho}$ is defined on the knot grid $\mathcal{T}^{\rho}$ with knot grid size $\rho$,
    the norms $\norm{\mathcal{I}}$ and $\norm{\mathcal{D}}$ are both related to the knot grid size $\rho$.
 Therefore, the sufficient condition for the consistency of the IGA-C method
    is followed.
 %%------------------------------------------------------------
% % Theorem
% %------------------------------------------------------------
% \begin{thm}[Sufficiency] \label{thm:sufficiency}
% Suppose $\mathcal{D}T = f \in C^0(\Omega_p)$ (Eq.~\pref{eq:bd_p}),
%        and $T_r \in \mathbb{S}_{\rho}$ (refer to Eq.~\pref{eq:rewrite-nurbs} and Remark 3) is the NURBS function
%        approximating the analytical solution $T$.
% If the interpolation operator $\mathcal{I}^\rho$~\pref{eq:inter_operator} and the differential operator
%    $\mathcal{D}$~\pref{eq:bd_p} are bounded when $\rho \rightarrow 0$,
%    then the IGA-C method applied on the boundary problem~\pref{eq:bd_p} is
%    consistency.
% \end{thm}

%------------------------------------------------------------
 % Theorem
 %------------------------------------------------------------
 \begin{lem}[Sufficiency] \label{thm:sufficiency}
 If the interpolation operator $\mathcal{I}^\rho$~\pref{eq:inter_operator}
    and differential operator $\mathcal{D}$~\pref{eq:bd_p} are both uniformly bounded when $\rho \rightarrow 0$,
    then the IGA-C method applied on the boundary problem~\pref{eq:bd_p} is
    consistency.
 \end{lem}

 Furthermore, the following lemma presents the necessary condition for
    the consistency of the IGA-C method.

 %---------------------------------------------------------------------------
 % Theorem
 %---------------------------------------------------------------------------

 \begin{lem}[Necessity] \label{thm:necessity}
     %Suppose $\mathcal{D}T = f \in C^0(\Omega_p)$ (Eq.~\pref{eq:bd_p}),
%        and $T_r \in \mathbb{S}_{\rho}$ (refer to Eq.~\pref{eq:rewrite-nurbs} and Remark 3) is the
%        NURBS function approximating the analytical solution $T$.
    If the IGA-C method applied on the boundary problem~\pref{eq:bd_p} is consistency,
    then the interpolation operator $\mathcal{I}^\rho$~\pref{eq:inter_operator} and the differential operator
    $\mathcal{D}$~\pref{eq:bd_p} are both uniformly bounded when $\rho \rightarrow 0$.
 \end{lem}

 \textbf{Proof:} We employ the method of proof by contradiction to
    show that $\mathcal{D}T_r$ is bounded when $\rho \rightarrow 0$.

 The consistency of the IGA-C method means that
 \begin{equation} \label{eq:consistency_cond}
        \mathcal{D}T_r \rightarrow \mathcal{D}T = f,\ \text{when}\ \rho \rightarrow
        0.
    \end{equation}
 By contradiction, suppose $\mathcal{D}T_r$ is not uniformly bounded
    when $\rho \rightarrow 0$, i.e.,
    $\norm{\mathcal{D}T_r}_{\mathbb{W}} \rightarrow \infty$,
    when $\rho \rightarrow 0$.
 Because $f$ is continuous,
    it is bounded on its domain $\Omega_p \cup \partial \Omega_p$.
 However, $\mathcal{D}T_r$ is unbounded when $\rho \rightarrow 0$.
 This violates the consistency condition~\pref{eq:consistency_cond}.
 So the hypothesis is not true,
        $\mathcal{D}T_r$ is uniformly bounded when $\rho \rightarrow 0$.
 That is, there exists a positive constant $C_r$ such that
 \begin{equation*}
    \norm{\mathcal{D} T_r}_\mathbb{W} \leq C_r ,
    \ \text{when}\ \rho \rightarrow 0.
 \end{equation*}

 Therefore, we have
 \begin{equation*}
    \norm{\mathcal{D}} = \sup_{\norm{T_r}_\mathbb{V} = 1}
    \{\norm{\mathcal{D}T_r}_\mathbb{W}\} \leq C_r,\ \text{when}\ \rho
    \rightarrow 0,
 \end{equation*}
 and (refer to Eq.~\pref{eq:inter_operator})
 \begin{equation*}
    \norm{\mathcal{I}^{\rho}} =
        \sup_{\norm{f}_{L^{\infty}}=1} \{\norm{\mathcal{I}^{\rho}f}_\mathbb{W}\}
        =
        \sup_{\norm{f}_{L^{\infty}}=1} \{\norm{\mathcal{D}T_r}_\mathbb{W}\}
        \leq C_r,\ \text{when}\ \rho \rightarrow 0.
 \end{equation*}
 It means that the interpolation operator $\mathcal{I}^\rho$~\pref{eq:inter_operator}
    and the differential operator
    $\mathcal{D}$~\pref{eq:bd_p} are both uniformly bounded when $\rho \rightarrow 0$.
 $\Box$

 \vspace{0.2cm}

 Based on Lemmas~\ref{thm:sufficiency} and~\ref{thm:necessity},
    the necessary and sufficient condition for the consistency of the IGA-C method is followed.

 %---------------------------------------------------------------------------
 % Theorem
 %---------------------------------------------------------------------------

 \begin{thm} [Necessity and Sufficiency] \label{thm:suf_nes}
  The IGA-C method applied on the boundary problem~\pref{eq:bd_p} is
    consistency,
    if and only if the interpolation operator $\mathcal{I}^\rho$~\pref{eq:inter_operator} and differential operator
    $\mathcal{D}$~\pref{eq:bd_p} are both uniformly bounded when $\rho \rightarrow 0$.
 \end{thm}

%-------------------------------------------------------------------------
% Section: Numerical Examples
%-------------------------------------------------------------------------

%-------------------------------------------------------------------------------------
%   Figure
%-------------------------------------------------------------------------------------
 \begin{figure}[!h]
    \centering
  \subfigure[Diagram of $\norm{T_r}_{L^{\infty}}$ v.s. $\ln{\rho_i}$ for the 1D problem.]{
    \label{subfig:1d_tr}
    \includegraphics[width = 0.45\textwidth]
        {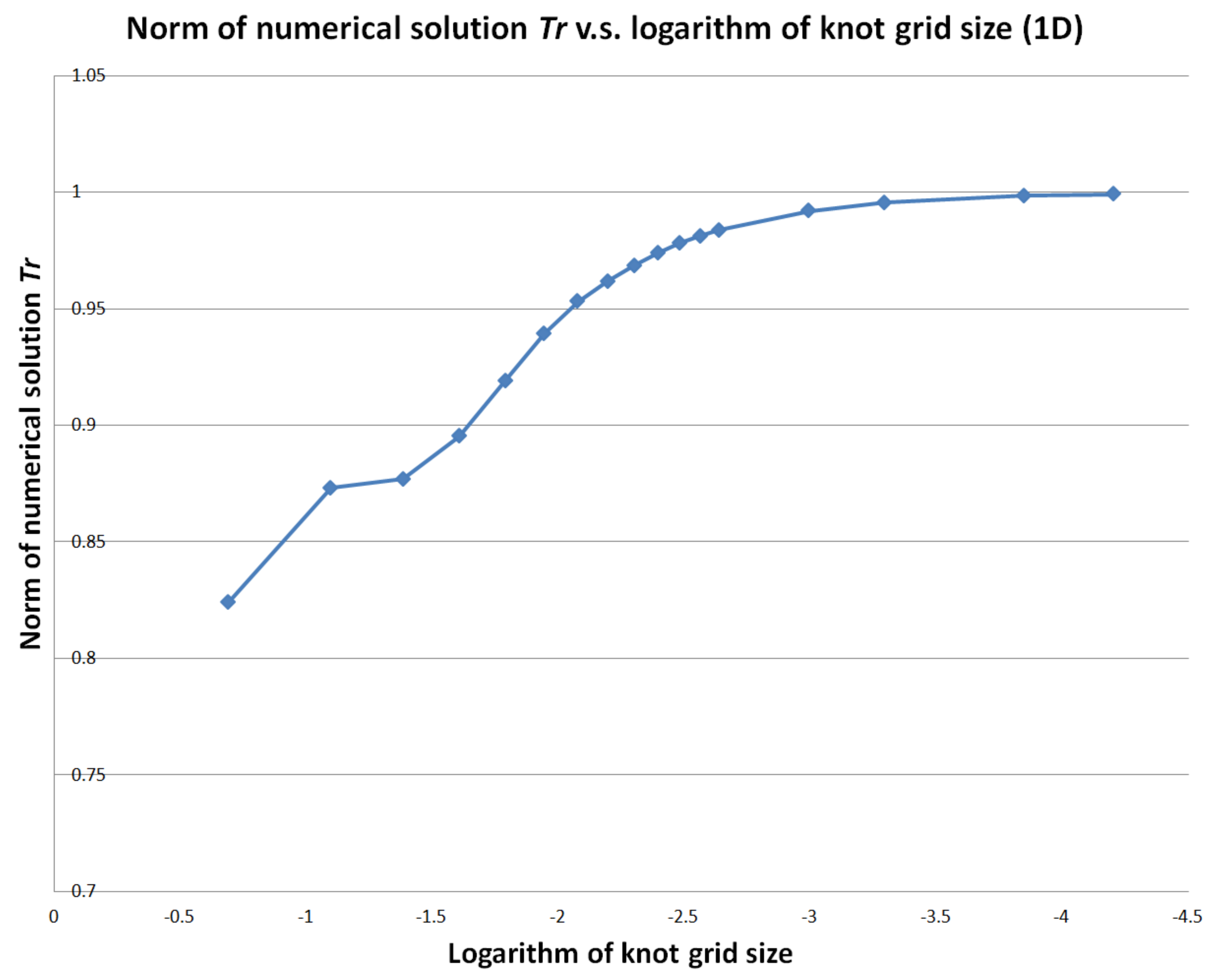}}
  \subfigure[Diagram of $\norm{\mathcal{D}T_r}_{L^{\infty}}$ v.s. $\ln{\rho_i}$ for the 1D problem.]{
    \label{subfig:1d_dtr}
    \includegraphics[width = 0.45\textwidth]
        {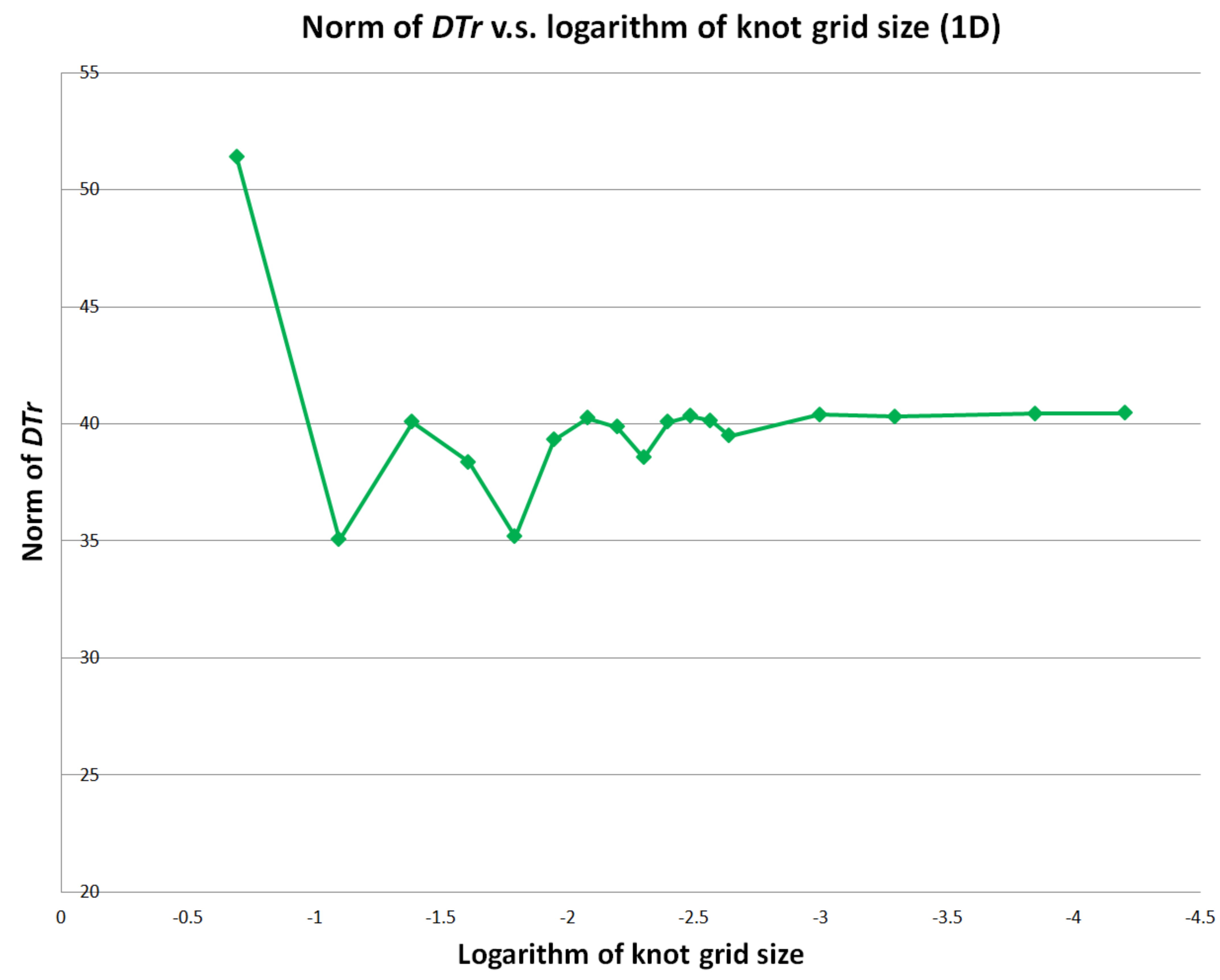}}
  \subfigure[Diagram of $\frac{\norm{\mathcal{D}T_r}_{L^{\infty}}}{\norm{T_r}_{L^{\infty}}}$ v.s.
                $\ln{\rho_i}$ for the 1D problem.]{
    \label{subfig:1d_dtr_tr}
    \includegraphics[width = 0.45\textwidth]
        {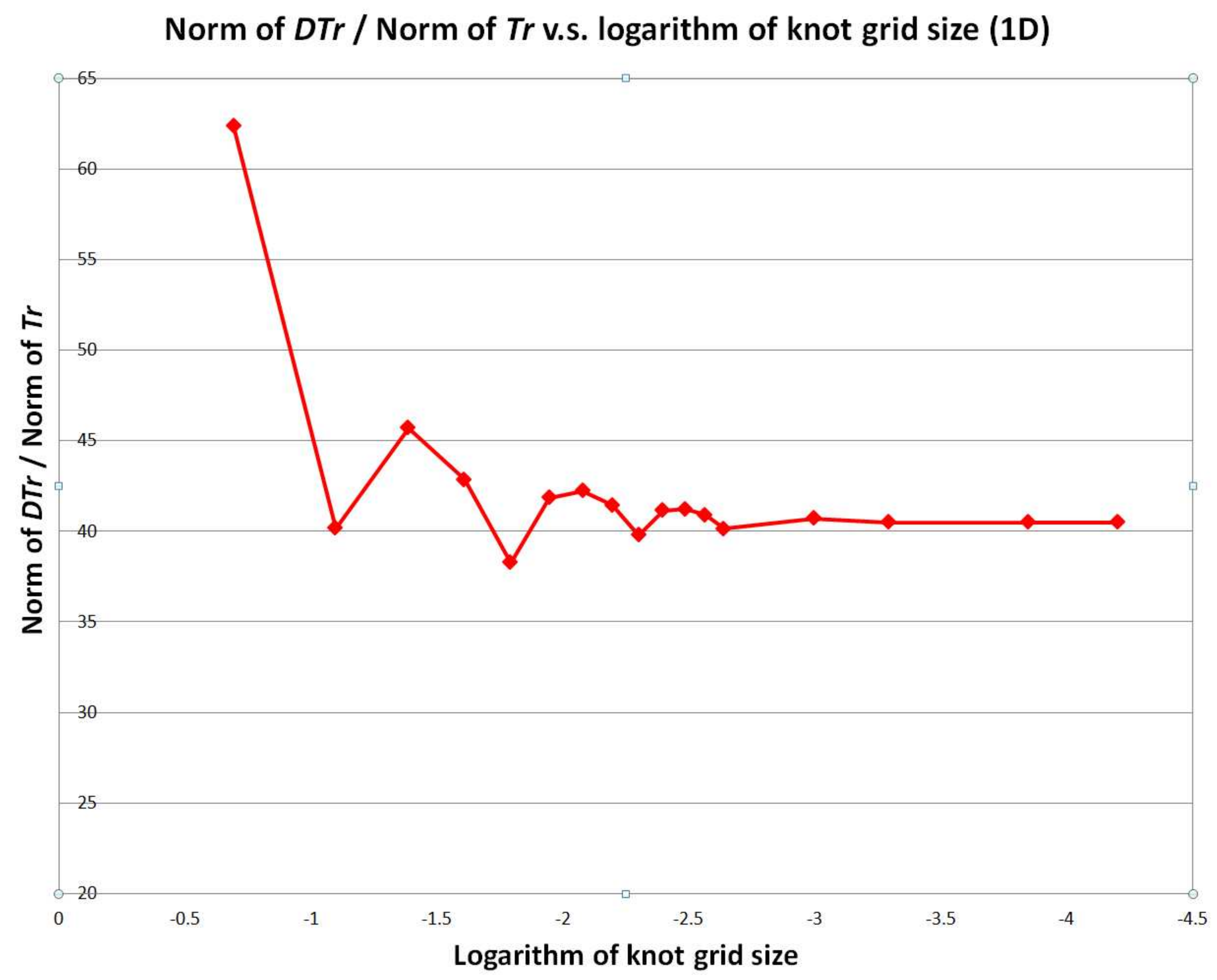}}
  \caption{In the case of one-dimensional source problem~\pref{eq:example_1},
  $\norm{T_r}_{L^{\infty}}, \norm{\mathcal{D}T_r}_{L^{\infty}}$, and the ratio
  $\frac{\norm{\mathcal{D}T_r}_{\infty}}{\norm{T_r}_{L^{\infty}}}$ are all uniformly bounded
  when the knot grid size sequence $\rho_k \rightarrow 0,\ (k \rightarrow \infty$).
}
  \label{fig:one_dim_example}
\end{figure}
%----------------------------------------------------------------------------------------------------------

\section{Numerical examples}
\label{sec:numerical_examples}

  In this section, some numerical examples are presented to illustrate the
    necessary and sufficient condition of the consistency of the IGA-C method.

 \textbf{Example 1: }
 Consider the following one-dimensional source problem:
 \begin{equation} \label{eq:example_1}
   \begin{cases}
   & -T'' + T = (1+4\pi^2) sin(2 \pi x), \quad x \in \Omega = [0,1],\\
   & T(0)=0,\ T(1) = 0.
   \end{cases}
 \end{equation}
 The analytical solution to the source problem is $T(x) = sin(2 \pi x)$.
 The physical domain is modeled by a cubic B-spline curve with control points
    $\{0, \frac{1}{3}, \frac{2}{3}, 1\}$ and knot vector $\{0\ 0\ 0\ 0\ 1\ 1\ 1\ 1\}$.
 So the initial knot grid size is $\rho_0 = 1$.
 To reduce the knot grid size,
    we uniformly insert $k, (k=1,2,\cdots)$ knots in $(0,1)$.
 And then, the knot grid size sequence is $\rho_k = \frac{1}{k+1}, k=0,1,2,\cdots$.

 In Fig.~\ref{fig:one_dim_example},
    three diagrams are demonstrated, that is,
    \begin{itemize}
      \item the norm of numerical solution $T_r$, i.e.,
            $\norm{T_r}_{L^{\infty}}$
            v.s. the logarithm of knot grid size, i.e., $\ln{\rho_k}$ (Fig.~\ref{subfig:1d_tr}),
      \item $\norm{\mathcal{D}T_r}_{L^{\infty}}$ v.s. $\ln{\rho_k}$
            (Fig.~\ref{subfig:1d_dtr}), and,
      \item $\frac{\norm{\mathcal{D}T_r}_{L^{\infty}}}{\norm{T_r}_{L^{\infty}}}$
            v.s. $\ln{\rho_k}$ (Fig.~\ref{subfig:1d_dtr_tr}).
    \end{itemize}
 It can be seen from the diagrams in Fig.~\ref{fig:one_dim_example} that,
    when $k \rightarrow \infty$ and $\rho_k \rightarrow 0$,
    the norm of the numerical solution $\norm{T_r}_{L^{\infty}}$ tends to the norm of the analytical solution, i.e., $\norm{T(x)}_{L^{\infty}} = \norm{sin(2 \pi x)}_{L^{\infty}} = 1$
    (Fig.~\ref{subfig:1d_tr}),
    and $\norm{\mathcal{D}T_r}_{L^{\infty}}$ tends to the norm of
    $f(x) = (1+4\pi^2) sin(2 \pi x)$~\pref{eq:example_1},
    i.e., $\norm{(1+4\pi^2) sin(2 \pi x)}_{L^{\infty}} = 1+4\pi^2$ (Fig.~\ref{subfig:1d_dtr}).
 Moreover, refer to Fig.~\ref{subfig:1d_dtr_tr},
    as an indicator of $\norm{\mathcal{D}}_{L^{\infty}}$,
    the ratio $\frac{\norm{\mathcal{D}T_r}_{\infty}}{\norm{T_r}_{L^{\infty}}}$ tends to
    $\frac{\norm{(1+4\pi^2) sin(2 \pi x)}_{L^{\infty}}}{\norm{sin(2 \pi x)}_{L^{\infty}}} = 1+4\pi^2$,
    when $\rho_k \rightarrow 0,\ (k \rightarrow \infty)$.
 Therefore, it is uniformly bounded as $\rho_k \rightarrow 0,\ (k \rightarrow \infty)$,
    which validates Theorem~\ref{thm:suf_nes}.

%-------------------------------------------------------------------------------------
%   Figure
%-------------------------------------------------------------------------------------
 \begin{figure}[!h]
    \centering
  \subfigure[Diagram of $\norm{T_r}_{L^{\infty}}$ v.s. $\ln{\rho_i}$ for the 2D problem.]
   {\label{subfig:2d_tr}
    \includegraphics[width = 0.45\textwidth]
        {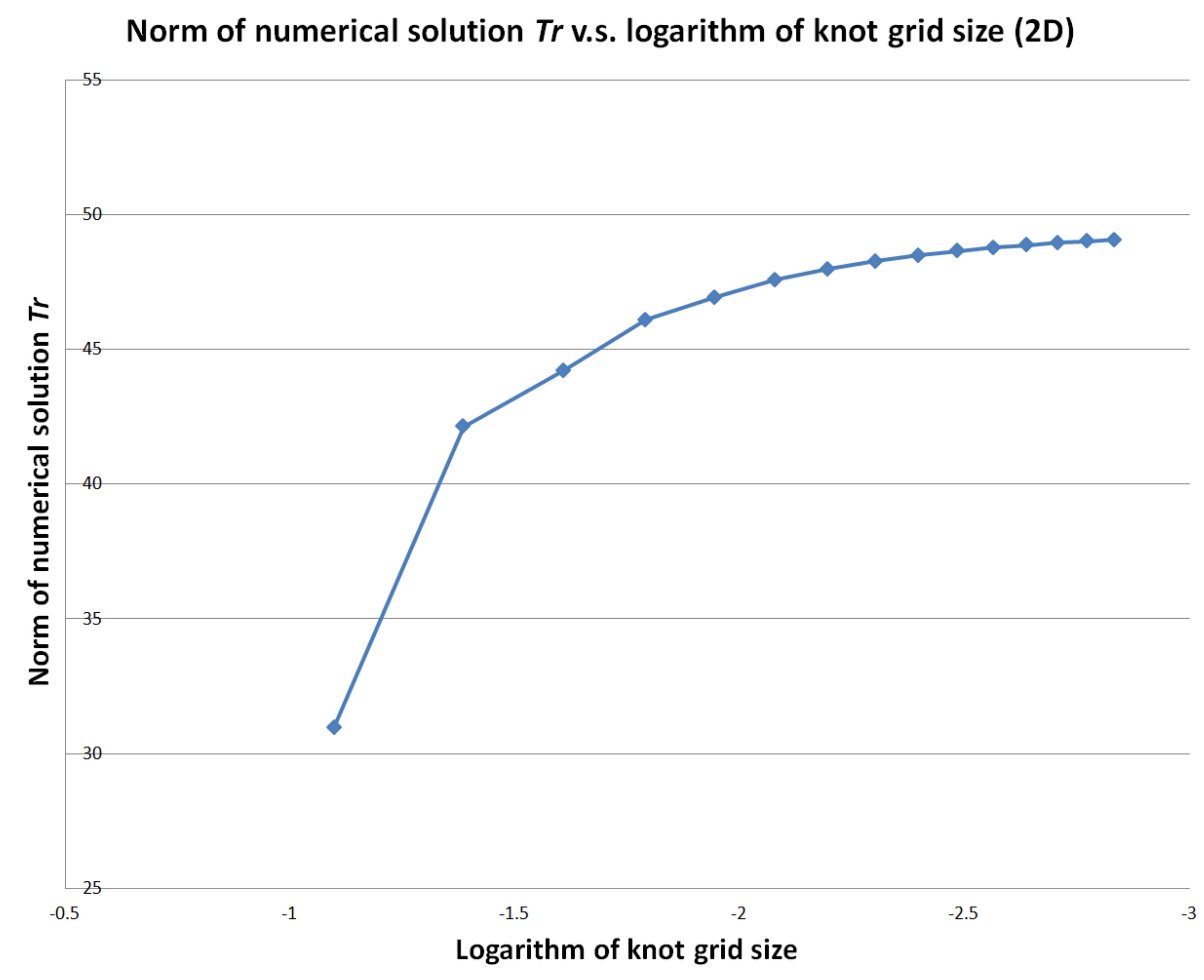}}
  \subfigure[Diagram of $\norm{\mathcal{D}T_r}_{L^{\infty}}$ v.s. $\ln{\rho_i}$ for the 2D problem.]
   {\label{subfig:2d_dtr}
    \includegraphics[width = 0.45\textwidth]
        {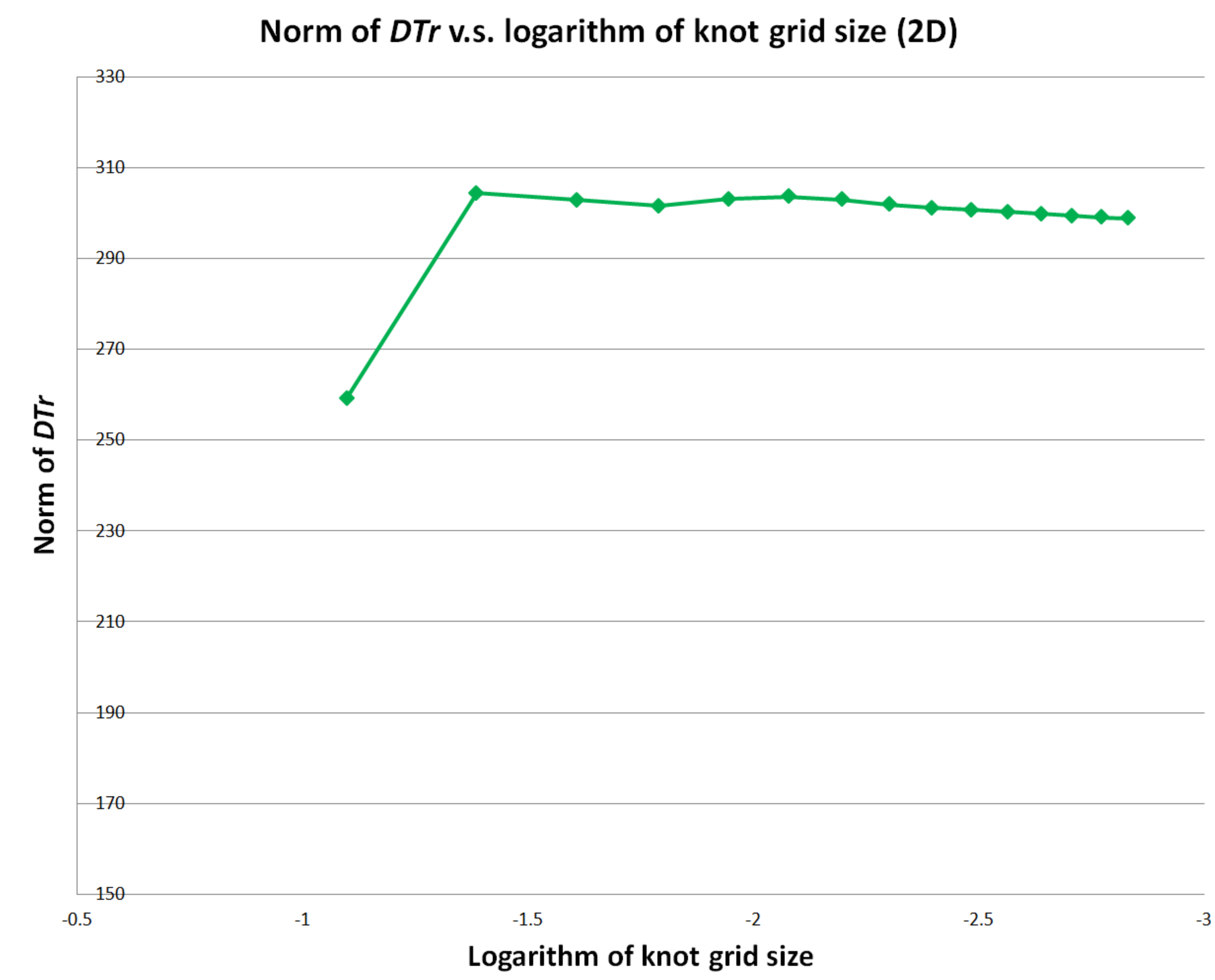}}
  \subfigure[Diagram of $\frac{\norm{\mathcal{D}T_r}_{L^{\infty}}}{\norm{T_r}_{L^{\infty}}}$ v.s.
                $\ln{\rho_i}$ for the 2D problem.]
   {\label{subfig:2d_dtr_tr}
    \includegraphics[width = 0.45\textwidth]
        {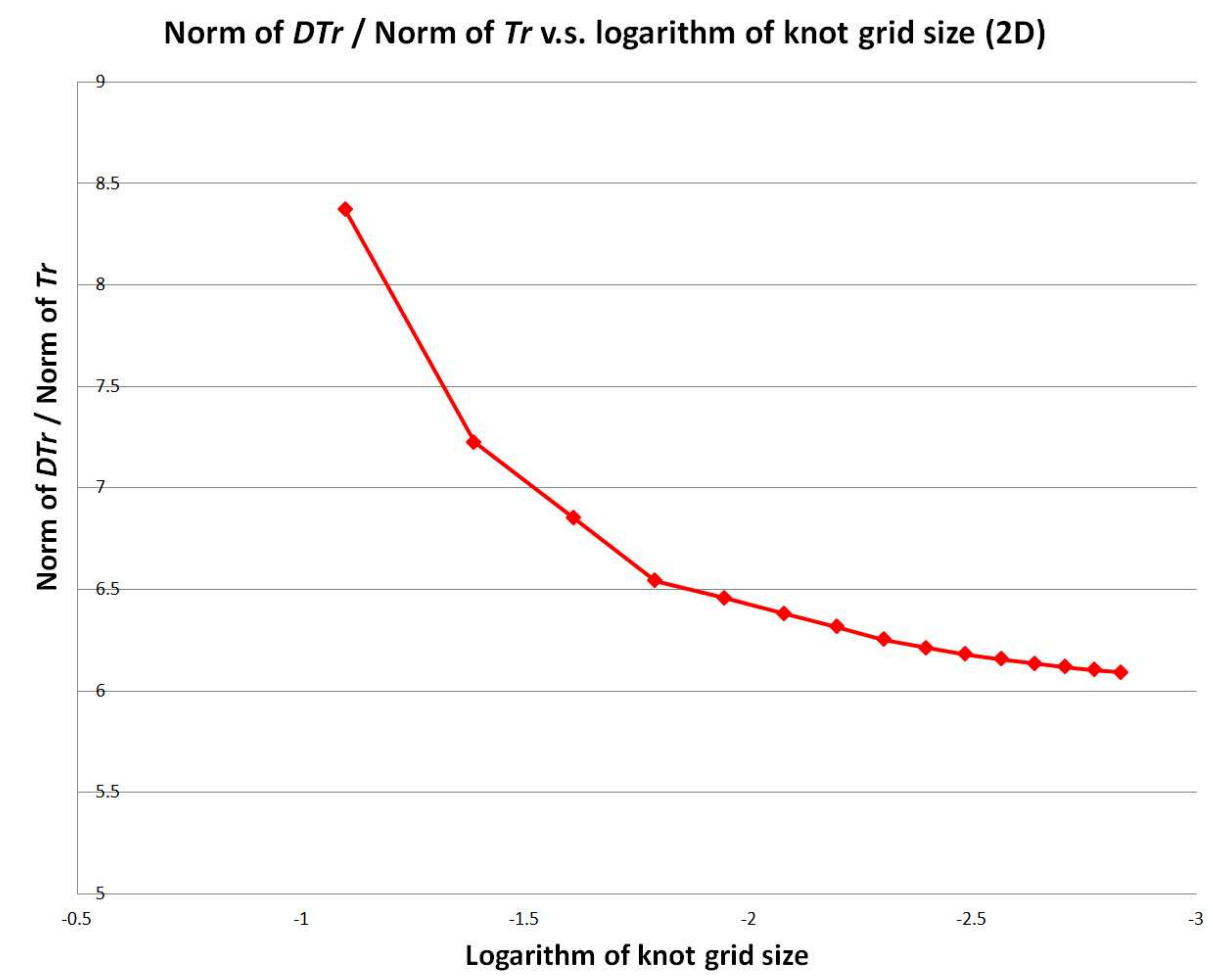}}
  \caption{In the case of two-dimensional source problem~\pref{eq:example_2},
  $\norm{T_r}_{L^{\infty}}, \norm{\mathcal{D}T_r}_{L^{\infty}}$, and the ratio
  $\frac{\norm{\mathcal{D}T_r}_{\infty}}{\norm{T_r}_{L^{\infty}}}$ are all uniformly bounded
  when the knot grid size sequence $\rho_k \rightarrow 0,\ (k \rightarrow \infty$).
}
  \label{fig:two_dim_example}
\end{figure}
%----------------------------------------------------------------------------------------------------------

 \textbf{Example 2: }
 The next example is a two-dimensional source problem:
     \begin{equation}
        \label{eq:example_2}
        \begin{cases}
        &-\Delta T + T = f,\  (x,y) \in \Omega \\
        &T|_{\partial {\Omega}} = 0,
        \end{cases}
    \end{equation}
    where,
    \begin{equation*}
        \begin{aligned}
        f = (3x^{4} - 67x^{2} - 67y^{2} + 3y^{4} + 6x^{2}y^{2} + 116)\sin(x)\sin(y) \\
        + (68x - 8x^{3} - 8xy^{2})\cos(x)\sin(y) \\
        + (68y - 8y^{3}-8yx^{2})\cos(y)\sin(x).
        \end{aligned}
    \end{equation*}
 And the analytical solution of the source problem~\pref{eq:example_2}
    is
    \begin{displaymath}
        T = (x^2 + y^2 - 1)(x^2 + y^2 - 16)\sin(x)\sin(y).
    \end{displaymath}

 The physical domain $\Omega$ in Eq.~\pref{eq:example_2} is a quarter of an annulus,
    which is represented by a cubic NURBS patch with $4 \times 4$ control points.
 The control points and weights of the cubic NURBS patch are listed in
    Tables~\ref{tbl:two-dim-ctrlpnt} and~\ref{tbl:two-dim-weights}, respectively.
 The knot vectors of the cubic NURBS patch along $u-$ and $v-$direction are, respectively,
    \begin{equation*}
        \begin{split}
            0\ 0\ 0\ 0\ 1\ 1\ 1\ 1, \\
            0\ 0\ 0\ 0\ 1\ 1\ 1\ 1.
        \end{split}
    \end{equation*}
 To make the knot grid size tend to $0$,
    we uniformly insert knots in the interval $(0,1)$ along $u-$ and $v-$directions, respectively.
 So, the knot grid sizes are $\rho_k = \frac{1}{k+1}, k=0,1,2,\cdots$.

 Fig.~\ref{fig:two_dim_example} shows the diagrams $\norm{T_r}_{L^{\infty}}$
    v.s. $\ln{\rho_k}$ (Fig.~\ref{subfig:2d_tr}),
    $\norm{\mathcal{D}T_r}_{L^{\infty}}$ v.s. $\ln{\rho_k}$ (Fig.~\ref{subfig:2d_dtr}),
    and $\frac{\norm{\mathcal{D}T_r}_{L^{\infty}}}{\norm{T_r}_{L^{\infty}}}$
    v.s. $\ln{\rho_k}$ (Fig.~\ref{subfig:2d_dtr_tr}) for the case of two-dimensional
    source problem~\pref{eq:example_2}.
 Similar as the case of one-dimensional problem,
    $\norm{T_r}_{L^{\infty}}$, $\norm{\mathcal{D}T_r}_{L^{\infty}}$,
    and $\frac{\norm{\mathcal{D}T_r}_{L^{\infty}}}{\norm{T_r}_{L^{\infty}}}$
    are all have limit when $\rho_k \rightarrow 0, (k \rightarrow \infty)$.
 So they are all uniformly bounded when $\rho_k \rightarrow 0, (k \rightarrow \infty)$.

%--------------------------------------------------------------------------
% Table
%--------------------------------------------------------------------------
 \begin{table}[!htb]
    \caption{Control points of the quarter of annulus}
    \label{tbl:two-dim-ctrlpnt} \centering
    \begin{tabular}{c  c  c  c  c}
    %\hline \multicolumn{2}{c}{Example I} &
    %\multicolumn{2}{c}{Example II} &
    %\multicolumn{2}{c}{Example III}\\
     \hline
    $i$ & $\bm{P}_{i,1}$ & $\bm{P}_{i,2}$ & $\bm{P}_{i,3}$ & $\bm{P}_{i,4}$ \\
     \hline
    1  & (1,0) & (2,0)  & (3,0) & (4,0)  \\

    2  & (1,2-$\sqrt{2}$) & (2, 4-2$\sqrt{2}$) & (3,6-3$\sqrt{2}$) & (4,8-4$\sqrt{2}$)\\

    3  & (2-$\sqrt{2}$,1) & (4-2$\sqrt{2}$,2)  & (6-3$\sqrt{2}$,3) & (8-4$\sqrt{2}$, 4)\\

    4  & (0,1) & (0,2) & (0,3) & (0,4)  \\
     \hline
    \end{tabular}
\end{table}

%---------------------------------------------------------------------------

%--------------------------------------------------------------------------
% Table
%--------------------------------------------------------------------------
\begin{table}[!htb]
    \caption{Weights for the quarter of annulus}
    \label{tbl:two-dim-weights} \centering
    \begin{tabular}{c  c  c  c  c}
    %\hline \multicolumn{2}{c}{Example I} &
    %\multicolumn{2}{c}{Example II} &
    %\multicolumn{2}{c}{Example III}\\
     \hline
    i & $w_{i,1}$ & $w_{i,2}$ & $w_{i,3}$ & $w_{i,4}$ \\
     \hline
    1  & 1 & 1  & 1 & 1  \\
    2  & $\frac{1+\sqrt{2}}{3}$ & $\frac{1+\sqrt{2}}{3}$ & $\frac{1+\sqrt{2}}{3}$ & $\frac{1+\sqrt{2}}{3}$ \\
    3  & $\frac{1+\sqrt{2}}{3}$ & $\frac{1+\sqrt{2}}{3}$  & $\frac{1+\sqrt{2}}{3}$ & $\frac{1+\sqrt{2}}{3}$ \\
    4  & 1 & 1 & 1 & 1  \\
     \hline
    \end{tabular}
\end{table}

%-------------------------------------------------------------------------------------
%   Figure
%-------------------------------------------------------------------------------------
 \begin{figure}[!h]
    \centering
  \subfigure[Diagram of $\norm{T_r}_{L^{\infty}}$ v.s. $\ln{\rho_i}$ for the 3D problem.]
   {\label{subfig:3d_tr}
    \includegraphics[width = 0.45\textwidth]
        {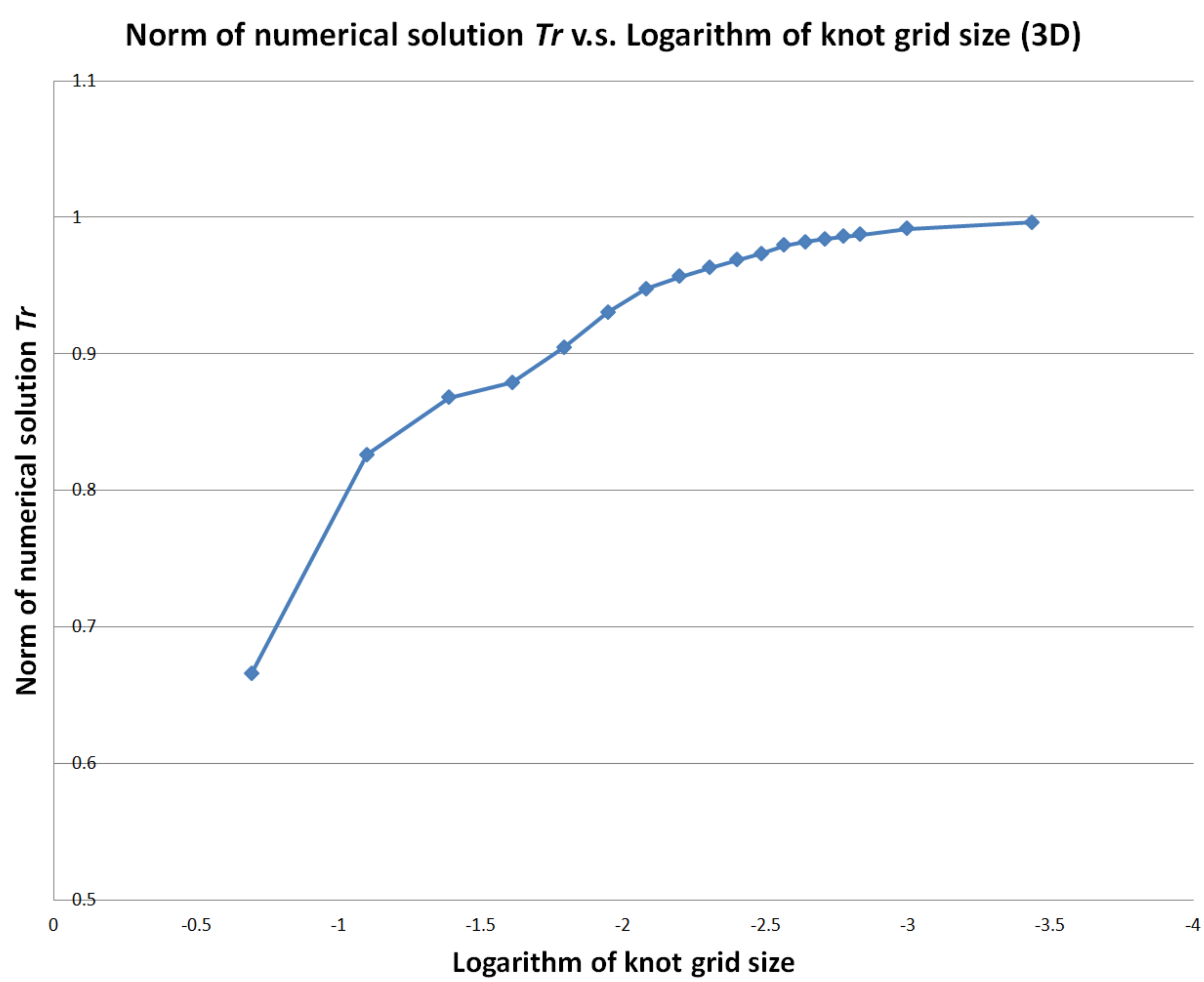}}
  \subfigure[Diagram of $\norm{\mathcal{D}T_r}_{L^{\infty}}$ v.s. $\ln{\rho_i}$ for the 3D problem.]
   {\label{subfig:3d_dtr}
    \includegraphics[width = 0.45\textwidth]
        {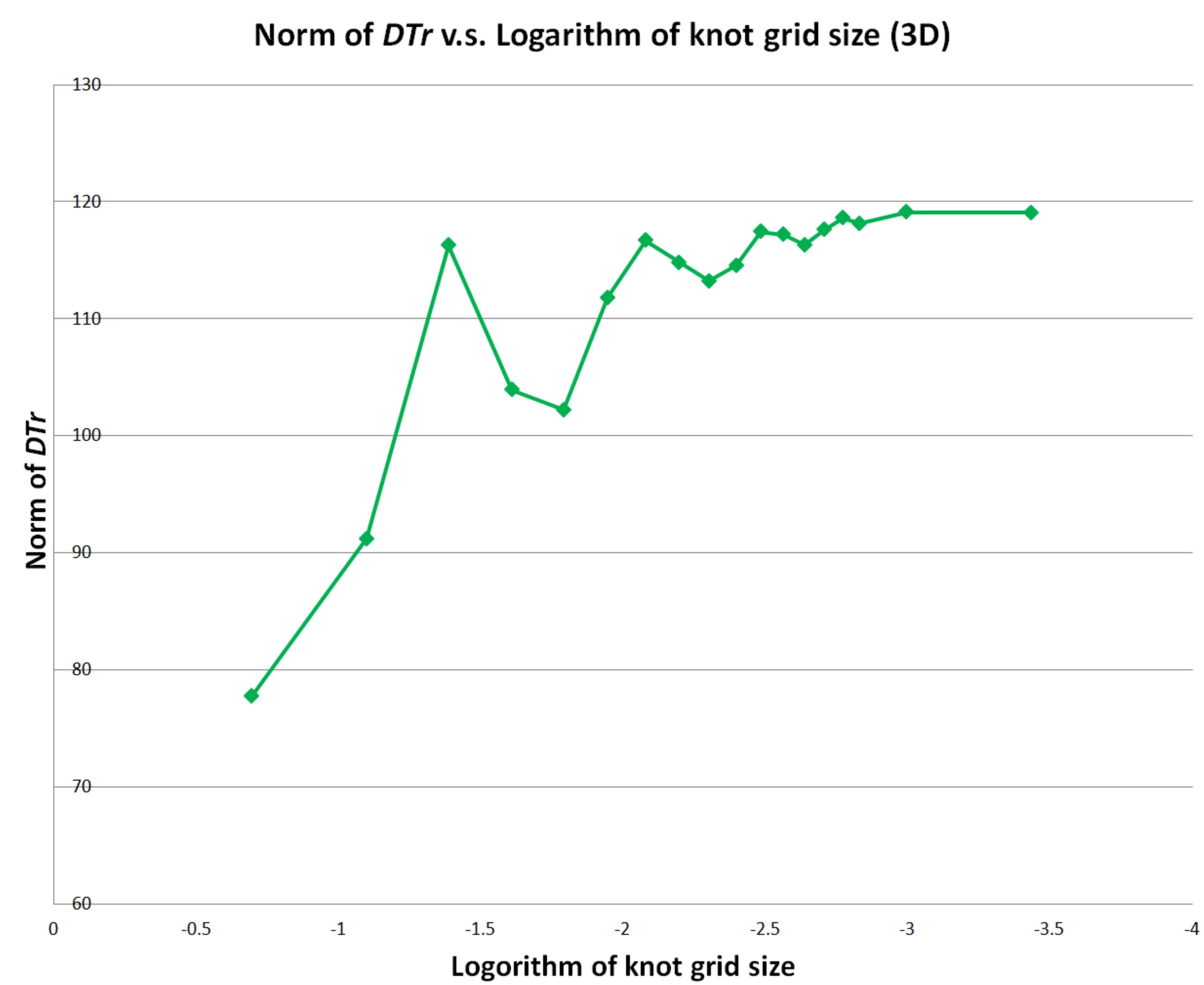}}
  \subfigure[Diagram of $\frac{\norm{\mathcal{D}T_r}_{L^{\infty}}}{\norm{T_r}_{L^{\infty}}}$ v.s.
                $\ln{\rho_i}$ for the 3D problem.]
   {\label{subfig:3d_dtr_tr}
    \includegraphics[width = 0.45\textwidth]
        {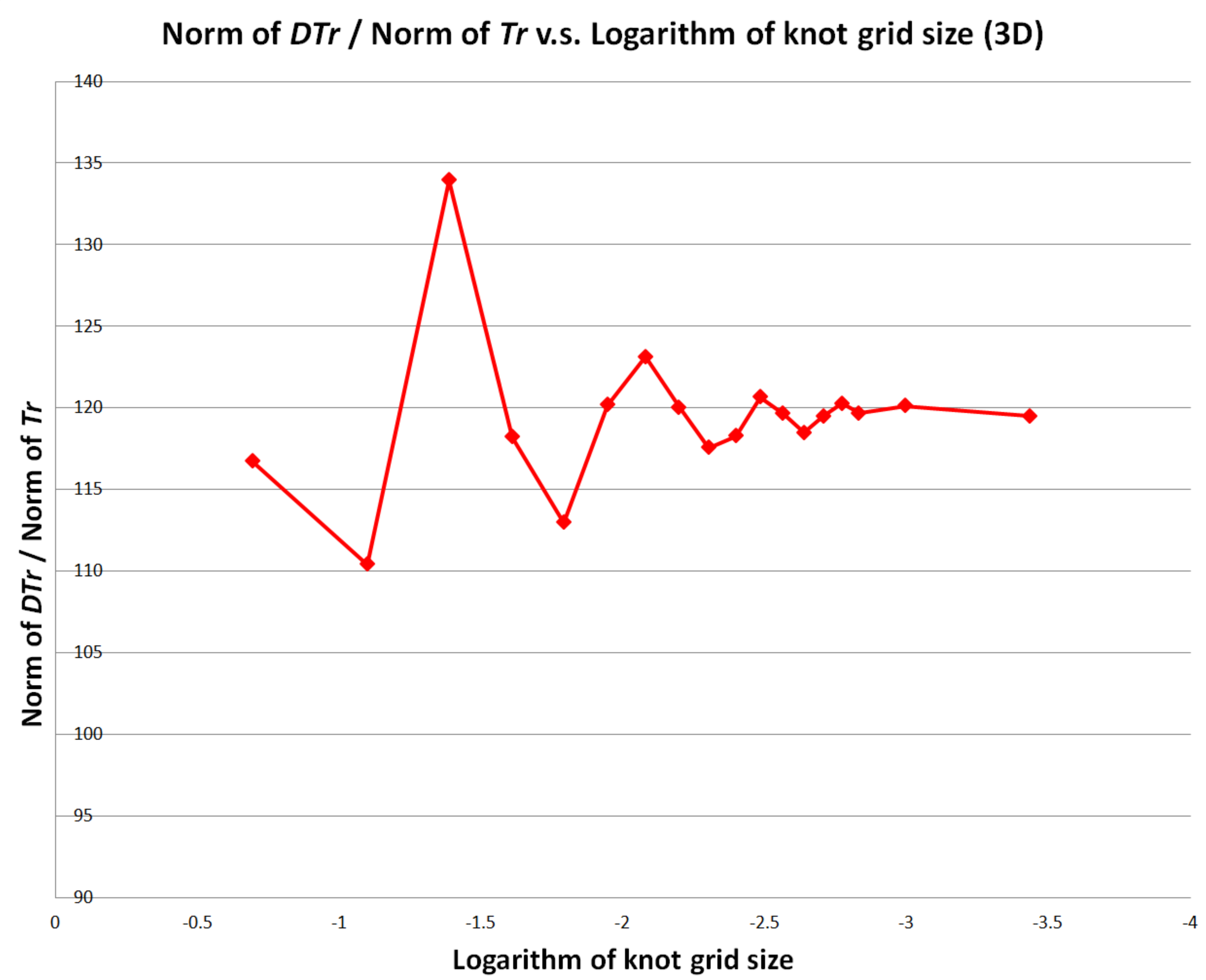}}
  \caption{In the case of three-dimensional source problem~\pref{eq:example_3},
  $\norm{T_r}_{L^{\infty}}, \norm{\mathcal{D}T_r}_{L^{\infty}}$, and the ratio
  $\frac{\norm{\mathcal{D}T_r}_{\infty}}{\norm{T_r}_{L^{\infty}}}$ are all uniformly bounded
  when the knot grid size sequence $\rho_k \rightarrow 0,\ (k \rightarrow \infty$).
}
  \label{fig:three_dim_example}
\end{figure}
%----------------------------------------------------------------------------------------------------------

 \textbf{Example 3:}
 The final example is a three-dimensional source problem:
 \begin{equation}
    \label{eq:example_3}
    \begin{cases}
    & -\Delta T + T = f,\  (x,y,z) \in \Omega, \\
    & T|_{\partial {\Omega}} = 0,
    \end{cases}
    \end{equation}
    where
    \begin{equation*}
    f = (1 + 12\pi^{2})\sin(2 \pi x)\sin(2 \pi y)\sin(2 \pi z),
    \end{equation*}
 and the analytical solution is,
    \begin{equation*}
    T = \sin(2 \pi x)\sin(2 \pi y)\sin(2 \pi z).
    \end{equation*}
 The physical domain $\Omega$ is modeled as a cubic trivariate B-spline solid
 with control points $\bm{P}_{ijk}=(\frac{i}{3}, \frac{j}{3}, \frac{k}{3}), i, j, k = 0,1,2,3$,
 and knot vectors along $u-$, $v-$, and $w-$directions, respectively,
 \begin{equation*}
        \begin{split}
            0\ 0\ 0\ 0\ 1\ 1\ 1\ 1, \\
            0\ 0\ 0\ 0\ 1\ 1\ 1\ 1, \\
            0\ 0\ 0\ 0\ 1\ 1\ 1\ 1.
        \end{split}
    \end{equation*}
 Similar as the one and two dimensional cases,
    the intervals $(0,1)$ along $u-$, $v-$, and $w-$directions are uniformly inserted knots, respectively.
 Therefore, the knot grid size $\rho_k = \frac{1}{k+1} \rightarrow 0,\ k = 0,1,2,\cdots$.

 The three diagrams, i.e., $\norm{T_r}_{L^{\infty}}$
    v.s. $\ln{\rho_k}$ (Fig.~\ref{subfig:3d_tr}),
    $\norm{\mathcal{D}T_r}_{L^{\infty}}$ v.s. $\ln{\rho_k}$ (Fig.~\ref{subfig:3d_dtr}),
    and $\frac{\norm{\mathcal{D}T_r}_{L^{\infty}}}{\norm{T_r}_{L^{\infty}}}$
    v.s. $\ln{\rho_k}$ (Fig.~\ref{subfig:3d_dtr_tr}) for the case of three-dimensional
    source problem~\pref{eq:example_3}
    are illustrated in Fig.~\ref{fig:three_dim_example}.
 From these diagrams, we can see that,
    when $\rho_k \rightarrow 0 (k \rightarrow +\infty)$,
    \begin{equation*}
      \norm{T_r}_{L^{\infty}} \rightarrow \norm{T}_{L^{\infty}} = 1, \quad
      \norm{\mathcal{D}T_r}_{L^{\infty}} \rightarrow \norm{f}_{L^{\infty}} = 1 + 12 \pi^2, \quad \text{and}, \quad
      \frac{\norm{\mathcal{D}T_r}_{L^{\infty}}}{\norm{T_r}_{L^{\infty}}} \rightarrow 1 + 12 \pi^2.
    \end{equation*}
 So they are all uniformly bounded when $\rho_k \rightarrow 0 (k \rightarrow \infty)$, too.

%-------------------------------------------------------------------------
% Section: Conclusion
%-------------------------------------------------------------------------
\section{Conclusions}
\label{sec:conclusion}

 In this paper, we developed the convergence order for the consistency and
    convergence of the IGA-C method,
    and then, deduced the necessary-and-sufficient condition for the
    consistency of the IGA-C method.
 Specifically, suppose $\mathcal{D}$ is the differential operator of
    a boundary value problem with $\mathcal{D}T=f$~\pref{eq:bd_p},
    a NURBS function $T_r$ is the numerical solution,
    and $\mathcal{I}^{\rho}$ is an interpolation operator such that
    $\mathcal{I}^{\rho} f = \mathcal{D}T_r$.
 First, the formula of the convergence order for the consistency of the IGA-C
    method is developed,
    which includes the norms of the operator $\mathcal{D}$ and $\mathcal{I}^{\rho}$.
 Then, the necessary-and-sufficient condition for the consistency of the
    IGA-C method is deduced.
 That is, the IGA-C method is consistency if and only if $\mathcal{D}$ and
    $\mathcal{I}^{\rho}$ are both uniformly bounded
    when $\rho \rightarrow 0$.
 These results will advance the numerical analysis of the IGA-C method.

%--------------------------------------------------------------------------------
% Section: Acknowledgement
%--------------------------------------------------------------------------------

\section*{Acknowledgement} This work is supported by the Natural Science Foundation of China
    (Nos. 61379072, 61202201).
    Dr. Qianqian Hu is also supported by the Open Project
    Program (No. A1305) of the State Key Lab of CAD\&CG, Zhejiang
    University.

% ----------------------------------------------------------------
%\bibliographystyle{plain}
%\bibliographystyle{unsrt}
%\bibliographystyle{abbrv}
%\bibliographystyle{alpha}
\bibliographystyle{elsart-num}
%\nocite{*}
\bibliography{isogeometric}

\end{document}